\documentclass[a4paper,11pt]{article}
\usepackage{amssymb}
\usepackage[english]{babel}
\usepackage{makeidx}
\usepackage[font=small]{caption}
\usepackage{amsmath}
\usepackage{amsthm}
\usepackage[margin=25mm]{geometry}
\usepackage{amsfonts}
\usepackage{amssymb}
\usepackage{mathrsfs}
\usepackage{graphicx}
\usepackage{enumerate}
\usepackage{indentfirst}
\usepackage{color}
\usepackage{authblk}
\makeindex
\definecolor{applegreen}{rgb}{0.55, 0.71, 0.0}

\newcommand{\M}{{\cal M}}

\newcommand{\intr}{\int_{\mathbb R^{N}}}
\newcommand{\ep}{\varepsilon}

\newcommand{\R}{\mathbb{R}}

\newcommand{\beq}{\begin{equation} }
\newcommand{\eqq}{\end{equation} }
\newcommand{\cuad}{{\sqcap\kern-.68em\sqcup}}

\newcommand{\rn}{\mathbb{R}^N}
\usepackage[utf8]{inputenc}

\renewcommand{\d}{\,\mathrm{d}}

\newtheorem{definition}{Definition}[section]
\newtheorem{teo}{Theorem}[section]

\newtheorem{proposition}{Proposition}[section]
\newtheorem{lemma}{Lemma}[section]
\newtheorem{corollary}{Corollary}[section]
\newtheorem{remark}{Remark}[section]
\newcommand{\bremark}{\begin{remark} \em}
\newcommand{\eremark}{\end{remark} }
\newtheorem{claim}{Claim}

\def\beeq{\begin{equation}}
\def\eeq{\end{equation}}
\newcommand{\begeqaet}{\begin{eqnarray*}}
\newcommand{\eneqaet}{\end{eqnarray*}}

\begin{document}
\begin{center}{\bf\Large 
Fundamental solutions and critical Lane-Emden exponents for nonlinear integral operators in cones
}\medskip
\bigskip
{\tiny
 
 {Gabrielle Nornberg}

Departamento de Ingeniería Matemática, Universidad de Chile

Beauchef 851, Santiago, Chile.

{\sl (gnornberg@dim.uchile.cl)}

\medskip
{Disson dos Prazeres} 

Departamento de Matemática, Universidade Federal de Sergipe

Av. Marechal Rondon s/n, Jd. Rosa Elze, São Cristóvão-SE, Brazil.

{\sl (disson@mat.ufs.br)}
\medskip

 {Alexander Quaas}

Departamento de  Matem\'atica,  Universidad T\'ecnica Federico Santa Mar\'{i}a

Casilla: V-110, Avda.\ Espa\~na 1680, Valpara\'{\i}so, Chile.

 {\sl (alexander.quaas@usm.cl)}
\medskip

 }
 
\end{center}

\medskip
\medskip
\medskip
\medskip

\begin{abstract}

In this article we study the fundamental solutions or ``$\alpha$-harmonic functions" for some nonlinear positive homogeneous nonlocal elliptic problems in conical domains, such as
\begin{eqnarray*}\label{ecbir1a1}
{\mathcal F }(u)=0\ \ \hbox{in} \ \ \mathcal{C}_\omega,\quad 
u=0\ \ \hbox{in} \ \ \rn\setminus \mathcal{C}_\omega ,\ \ 
\end{eqnarray*}
where $\omega$ is a proper $C^2$ domain in $S^{N-1}$ for $ N\geq 2$, $\mathcal{C}_\omega:=\{x\,:\,x\neq 0, {|x|^{-1}}x\in \omega\}$ is the
cone-like domain related to $\omega$, and ${\mathcal F }$ is an extremal fully nonlinear integral operator. 
We prove the existence of two fundamental solutions that are homogeneous and do not change signs in the cone; one is bounded at the origin and the other at infinity.

As an application, we use the fundamental solutions obtained to prove Liouville type theorems in cones for supersolutions of the Lane-Emden-Fowler equation in the form
\begin{eqnarray*}\label{eq 0.2}
{\mathcal F }(u)+u^p = 0\ \ \hbox{in} \ \ \mathcal{C}_\omega, \quad 
u=0\ \ \hbox{in} \ \  \rn\setminus \mathcal{C}_\omega.\ \ 
\end{eqnarray*}
We also prove a generalized Hopf type lemma in domains with corners. 
Most of our results are new even when ${\mathcal F }$ is the fractional Laplacian operator. 

 \end{abstract}
\date{}

\bigskip

{\textbf{Keywords.} Fully nonlinear equations, nonlinear integral operators; fundamental solutions; Liouville type theorems.\medskip
	
{\textbf{MSC2020.} 47G20, 45K05, 35B53, 35J15, 35J60.

\setcounter{equation}{0}
\section{ Introduction}

Fundamental solutions are one of the most basic tools in PDEs, by being particularly useful in fully nonlinear scenarios.
Build from a convenient symmetry, these solutions play a role as optimal barriers to achieve sharp results in diverse complex settings, some of them discussed below.  Our goal in this paper is twofold: first constructing fundamental solutions, and next applying them to establish nonlinear Liouville theorems for Lane-Emden type nonlocal equations defined in cones. 

Our analysis starts with the study of fundamental solutions, or $\alpha$-harmonic functions, for some nonlinear nonlocal elliptic problems defined in cone-like domains of the form
\begin{equation}\label{ecbir1a1}
 \left\{
\begin{array}{rllll}
{\mathcal F }(u)&=&0& \rm { in} & \;\mathcal{C}_\omega,\\
 u&=&0& \rm {in} & (\mathcal{C}_\omega)^c ,
\end{array}
\right.
\end{equation}
where $\omega \subsetneq S^{N-1}$ is a $C^2$ domain for $N\geq 2$, $\mathcal{C}_\omega:=\{x\, : \,x\neq 0, {|x|^{-1}}x\in \omega\}$ is the cone corresponding to $\omega$. Here, ${\mathcal F }$ is an extremal fully nonlinear integral operator such as
\begin{align}\label{Pucci}
{\cal M}^+(u)=\sup_{K}L_K(u) \quad\mbox{or}\quad
{\cal M}^-(u)=\inf_{K}L_K(u),
\end{align}
where the supremum and infimum are taken under positive even kernels $K:\R^N\to\R$ satisfying
$
\frac{\lambda}{|y|^{N+2\alpha}}\le K(y)\le \frac{\Lambda}{|y|^{N+2\alpha}}
$
for $\Lambda\ge \lambda>0$, $\alpha\in (0,1)$, and $L_K(u)$ is a linear operator as
\begin{align}\label{LK}
	L_K(u)(x)=\int_{\R^N}\{u(x+y)+u(x-y)-2u(x) \} K(y)\d y.
\end{align}

By fundamental solutions we mean two homogeneous viscosity solutions of \eqref{ecbir1a1} that do not change sign in the cone, such that one of them is bounded at the origin whereas the other is bounded at infinity.
To the best of our knowledge, the applications of fundamental solutions presented in this study are new,  even in the case ${\mathcal F }= -(-\Delta)^{\alpha}$, where  $(-\Delta)^{\alpha}$ denotes the fractional Laplacian:
 \begin{equation*}\label{fracLapl}
(-\Delta)^\alpha u(x) = C_{N,\alpha}
\int_{\R^N}\{u(x+y)+u(x-y)-2u(x) \} \frac{ 1}{|y|^{N+2\alpha}} \d y,
\end{equation*} 
for a constant $C_{N,\alpha}$ depending on $N$ and $\alpha$ that we omit for simplicity.
This is the most basic nonlocal operator we may consider, followed by linear operators as in \eqref{LK} with $-(N+2\alpha)$ homogeneous kernels. But we can go further allowing $\mathcal{F}$ to be of Isaacs type, by then including a large class of nonlinear integral operators, see Section \ref{section preliminaries}.

In a nutshell, nonlocal operators appear as  infinitesimal generators of stochastic Levy processes and optimal controls \cite{B, O, S}. Some references of nonlocal diffusion phenomena include \cite{La} to Particle Models in Physics, \cite{CR} to Nonlinear Reaction-Diffusion for Population Biology, \cite{val} for a big list of applications of the so called anomalous diffusion, just to quote a few.

In the nonlinear scenario, countless contributions have been made since the seminal series of papers \cite{CS1, CS2, CS3} by Caffarelli and Silvestre, where they introduced and developed the basis of the regularity theory for fully nonlinear integral equations. Fundamental solutions and Liouville type applications in the spirit we consider in this work were established in  \cite{FQ, FQDCDS} for nonlinear integral equations defined in $\mathbb{R}^N$, respectively for Pucci and Isaacs operators. 

It is important to recall some closely connected literature involving the local case in a fully nonlinear perspective. Labutin  \cite{labu} and Cutri-Leoni \cite{CL}  studied the fundamental solutions for Pucci operators in the whole $\mathbb{R}^N$. In these studies, the motivation is centered on the removal of singularities, Hadamard type three-sphere properties, and Liouville type nonexistence theorems.
In a more general setting, existence of fundamental solutions  in the whole $\mathbb{R}^N$ was established by Armstrong, Sirakov and Smart in \cite{ASS1}, together with applications in stochastic differential games.
As far as cone-like domains are concerned, in \cite{ASS} the same authors established existence of fundamental solutions for a general class of fully nonlinear elliptic operators, by extending previous results due to Miller \cite{Mi}.  In \cite{ASS}, applications of these fundamental solutions in cones imply Phragmen-Lindel\"of theorems, Hopf type lemma and Picard-Bouligand type principles.

In what refers to the nonlocal case, fundamental solutions for the fractional Laplacian, so called $\alpha$-stable or $\alpha$-harmonic, were obtained for cone like-domains in
\cite{bob3, bob1, bob2}, see also  \cite{proba} to the probability point of view and references therein.
Asymptotic results of these $\alpha$-harmonic functions in cones as $\alpha \to 1$ were studied by Terracini, Tortone and Vita in \cite{terra}.
More recently, Fern\'andez-Real and Ros-Oton used this approach in cones to study a nonlocal thin one-phase free boundary problem, see \cite{Ros}.

\medskip

In the sequel we enunciate our first main result regarding the existence of fundamental solutions for Pucci nonlocal operators in conical domains. It also holds for Isaacs type operators as in \eqref{Isaacs}, see next section for the precise definition.
The solutions we obtain do not change sign; moreover, they are homogeneous and bounded either at the origin or at infinity. Solutions are understood in the viscosity sense and are continuous in the (open) cone as in Definition~\ref{DefVisc} ahead. We point out that the fundamental solutions in $\R^N$ play a role in the assumptions. 

From the results of ~\cite{FQ} it is known that there exist \textsl{dimension-like numbers} $\tilde N^{\pm} > 0$ (depending on the ellipticity constants) such that the function $|x|^{2\alpha - \tilde N^{\pm}}$ is a fundamental solution associated to $\M^{\pm}$. 
Notice that the inequalities $\tilde N^+ > 2\alpha$ or  $\tilde N^+ < 2\alpha$ depend on the ellipticity constants, essentially, $\tilde N^+ >2\alpha$ holds  if the ellipticity constants are not too far away, whereas the opposite is true otherwise. In addition, $\tilde N^+=2\alpha$ is responsible for producing logarithmic solutions for $\M^+$ as in \cite{FQ}, see also \cite{FQDCDS}. Recall $ N\geq2$, so  $\tilde N^-\geq N>2\alpha$. 

In this article we will focus mostly on the case $\tilde N^+>2\alpha$, which is also 
a natural assumption in practical applications.  Nevertheless, we discuss all cases in the next theorem for a given $\omega\subset S^{N-1}$ fixed proper and $C^2$ domain in $S^{N-1}$.

\begin{teo}\label{tm 0.1} Let $\alpha\in (0,1)$ and $\mathcal{F}$ be as in \eqref{Pucci}--\eqref{LK}.
\begin{enumerate}[(a)]
	\item Assume $\tilde N^+>2\alpha$. Then there exists a unique constant $0<\beta^+(\omega, \mathcal{F})<N$ such that problem \eqref{ecbir1a1} possesses a fundamental positive solution  $\phi^+$ which is $-\beta^+(\omega, \mathcal{F})$ homogeneous.

	\item  If either $\tilde N^+< 2\alpha$ or $\omega\subset S^{N-1}_+$, there exists a unique constant $-2\alpha<\beta^-(\omega, \mathcal{F})<0$ so that \eqref{ecbir1a1} possesses a fundamental positive solution  $\phi^-$ which is $-\beta^-(\omega, \mathcal{F})$ homogeneous.
\end{enumerate}	
\end{teo}
We emphasize that the most interesting applications rely on $\omega\subset S^{N-1}_+$ and $\tilde N^+>2\alpha$, as shown below, where both fundamental solutions exist for $\mathcal{F}$. The particular case of the Pucci operator $\M^-$ does not require any hypothesis for the existence of $\beta^+(\omega,\mathcal{M}^-)$, cf.\ Lemma \ref{lemma well defined}.

The proof of Theorem \ref{tm 0.1} is inspired by ideas from \cite{ASS}, that is, the use of degree theory in an approximated weighted eigenvalue problem that blows up. Subsequently, a normalized function converges  to the fundamental solution as the parameter becomes extremal. These ideas also appear in some versions of Krein-Rutman theorem.
This  approximated eigenvalue problem is far from evident  in the nonlocal case. For example, in the auxiliary problem used in Lemma~\ref{le 2.2}, the nontrivial  $g$ function plays a delicate role.  Furthermore, crucial barriers require global upper (or lower) estimates in all complements of the domain to use comparison in that domain, which brings many technical difficulties compared with the local case. 

Now, we discuss an application to a fully nonlinear nonlocal Lane-Emden equation type in cones.  We are interested in studying existence and nonexistence results of the following nonlocal elliptic problem concerning a continuous function $u$ in the cone $\mathcal{C}_\omega$:
\begin{equation}\label{eq 0.002}
 \left\{
\begin{array}{rllll}
{\mathcal F }(u)+u^p&\leq &0 & \rm { in} &\; \mathcal{C}_\omega ,\\
 u&\geq &0& \rm {in} & (\mathcal{C}_\omega)^c .
\end{array}
\right.
\end{equation}

Before we present our contributions to equation \eqref{eq 0.002}, we discuss some previous results for the local case. Two pioneering works go back to Bandle, Levine \cite{ban2} and Bandle, Esssen  \cite{ban2} who established the first results of this type. Thereafter, Berestycki, Capuzzo-Dolcetta, and Nirenberg \cite{BCN} established sharp results for any cone (even without regularity).

Some extensions can be found in Section 5 of  \cite{AS}, as well as the references therein. In the special case of the half space, Leoni  in \cite{Leoni} proved  a precise estimate of the homogeneity of the fundamental solutions by computing the Pucci operator for special functions. Meanwhile, these precise estimates are a widely open problem in the nonlocal case, except for the case of the fractional Laplacian. 

Now, as a first consequence of the fundamental solutions, we derive the following  Liouville theorems for Lane-Emden equations driven by fully nonlinear nonlocal operators in cones.  This is the first result in this setting, as far as we know, for nonlocal operators even in the simplest case of the fractional Laplacian. To be coherent with our hypotheses and statements, we split our theorems into positive and nonpositive powers of $p$.
\begin{teo}\label{tm 0.2}
Assume $\tilde N_+ > 2\alpha$ and let $\beta^+=\beta^+(\omega,\mathcal{F})>0$ be the constant from Theorem \ref{tm 0.1}(a). If 
	 $0< p\leq \frac{\beta^++2\alpha}{\beta^+}$ 
then problem \eqref{eq 0.002}
has no positive solutions in the cone $\mathcal{C}_\omega$.
\end{teo}

Let us now fix the operator $\mathcal{F}=\M^-$ as in \cite{Leoni} and consider the supersolutions problem
\begin{equation}\label{eq 0.002M-}
	\mathcal{M}^-(u)+u^p\le 0 \quad \textrm{in } \mathcal{C}_\omega \, , \quad u\ge 0\;\; \textrm{ in } \mathbb{R}^N\setminus\mathcal{C}_\omega \, .
\end{equation}

\begin{teo}\label{tm 0.3} Assume $\tilde N^+< 2\alpha$ or $\omega\subset S^{N-1}_+$ and let $\beta^-=\beta^-(\omega,\mathcal{M}^-)<0$ be the constant from Theorem \ref{tm 0.1}(b). If 
$\frac{\beta^-+2\alpha}{\beta^-}\le p \le 0$ 
then problem \eqref{eq 0.002M-}	has no positive solutions in the cone $\mathcal{C}_\omega$.
	Moreover, any supersolution of \eqref{eq 0.002M-} for $p<\frac{\beta^-+2\alpha}{\beta^-}$ is unbounded at infinity. 
\end{teo}

In the case of the fractional Laplacian, as far as the equality signs in \eqref{eq 0.002} are concerned, we mention two different approaches:  in \cite{Xia} some results for the half space were developed through monotonicity arguments; while in \cite{FW} the Caffarelli-Silvestre extension was employed to starshaped like unbounded domains including $\mathbb{R}^N_+$.
However, ours seems to be the first result obtained for supersolutions, for which improved exponents are expected.
Moreover, we take into account the structure of the general cone through its respective fundamental solution.

It noteworthy that nonlinear Liouville type theorems for Lane-Emden equations in the special case of the half space are neatly connected with a priori bounds for nonconvex bounded domains, a feature considerably unknown in the case of a general fully nonlinear operator. 
Given their nonvariational nature, topological methods are typically employed to reach existence, hence the importance of deriving a priori estimates.
A common strategy to address this challenge relies on the well-known blow-up argument, which in turn requires a Liouville type theorem in $\mathbb R^N$ or $\mathbb R_+^N$ after blowing up; see \cite{GS} for more details and \cite{Revista} for nonlocal counterparts. 

\smallskip

In the model case of the fractional Laplacian operator defined in the half space $\mathcal{C}_\omega=\mathbb{R}_+^{N}$ for $\omega=S^{N-1}_+$, we can use the Kelvin transform\footnote{We recall that the Kelvin transform of $u$ is defined as $\bar{u}(x)=|x|^{-N+2\alpha}u(\bar{x})$, where $\bar{x}=x|x|^{-2}$. It satisfies
		$\Delta^{\alpha} \bar{u} (x)=|x|^{-N-2\alpha} \Delta^\alpha u (\bar{x})$ by using the identity $|\bar{x}-\bar{y}| \, |x|\,|y|=|x-y|$.
} 
and the fact that $(x_n)_+^\alpha$ is $\alpha$-harmonic in the half space to find that $(x_n)_+^\alpha/|x|^N$ is also $\alpha$-harmonic in $\mathbb{R}_+^{N}$, therefore
$\beta^+=N-\alpha$ and $\beta^-=-\alpha$.
\smallskip

As a corollary of our Liouville results in the case of the fractional Laplacian and the discussion above, we obtain an explicit characterization:
\begin{corollary} Assume $N\ge 2$.
For any real $p$ so that $-1\le p\leq \frac{N+\alpha}{N-\alpha}$, problem
\begin{align}
 -(-\Delta)^{\alpha}u+u^p\leq 0 \textrm{ in }  \mathbb R_+^N , \quad u\geq 0 \textrm{ in } \mathbb R^N_-,
\end{align}
possesses no positive solutions.
\end{corollary}

The critical exponent $\frac{N+\alpha}{N-\alpha}$ appearing above is new, up to our knowledge, in any context involving the  fractional Laplacian. 
When $\alpha=1$ this exponent corresponds to Brezis-Turner exponent that is $\frac{N+1}{N-1}$, see \cite{BT}.
On the other hand, the negative critical exponent $-1$ coincides with the one found in the more recent studies regarding the sublinear problem \cite{Leoni, moroz}. 
In what concerns the subcritical negative scenario, for any $p<-1$, it follows from the calculations in \cite{val, RSduke} that the function $(x_N)_+^{\alpha \gamma}$ for suitable $\gamma$ is an explicit positive supersolution in the half space.

As a second application, we explore the genuine shape of the second  fundamental solution 
$\phi^-$ near the origin  in order to obtain a generalized Hopf lemma at corners. 
We recall that $\phi^-$ is $-\beta^->0$ homogeneous, in particular bounded at the origin with $\phi^-(0)=0$. Thus, the statement will say that any bounded positive supersolution having a minimum $0$ at the vertex of the cone must be nontangential at this point. Notice that if the cone is strictly included in the half space then no interior tangent ball at zero completely contained in the cone is admissible. 
\begin{teo}\label{Hopf 0}
Let $\mathcal{F}$ as in \eqref{Pucci}--\eqref{LK} and $u$ be a bounded solution of 
\begin{align}
\mathcal{F}(u)\leq 0 \;\;\textrm{ and } \;\; u> 0 \textrm{ \,in }\,\mathcal{C}_\omega, \quad u=0 \textrm{ \,in }\,(\mathcal{C}_\omega)^c, \,\,u(0)=0.
\end{align}
Then, for any $\omega_0\subset\subset \omega\subset S^{N-1}_+$, there exists a constant $C>0$ depending only on $\lambda,\Lambda,N$, $\omega_0$, and $\mathrm{dist}(\omega_0,\partial\omega)$, such that, for any $\beta \in (-2\alpha, \beta^-(\omega, \mathcal{F}))$,
\begin{align}\label{est Hopf origem}
u(te)\geq Ct^{-\beta} \textrm{ \, as $t\rightarrow 0$, \; for all $e\in \omega_0$}.
\end{align}
\end{teo}

We notice that the method we use permits to treat weighted equations up to a natural change in the exponent; here we avoid the weight for simplicity matters.

\medskip

The remainder of this paper is organized as follows. Section \ref{section preliminaries} is devoted to some preliminary results concerning our extremal integral operators. In Section \ref{section fundamental} we prove the existence of fundamental solutions. Finally, we use the fundamental solutions obtained to prove our Liouville theorem in Section \ref{section Liouville}, together with an asymptotic behavior analysis.

 \setcounter{equation}{0}
  \section{Preliminaries}\label{section preliminaries}
  
  We start this section by defining the class ${\cal E}$ of integral operators we will consider throughout the text.
  Let $K:\R^N\to\R$ be a positive even function satisfying
  \begin{align}\label{K}
 \frac{\lambda}{|y|^{N+2\alpha}}\le K(y)\le \frac{\Lambda}{|y|^{N+2\alpha}},
  \end{align}
where $N\ge 2$, $\Lambda\ge \lambda>0$ and $\alpha\in (0,1)$. 
 For such a $K$ and for a suitable function $u$,  it is useful to denote the linear operator $L_K(u)$ in \eqref{LK}  as
  \begin{align}\label{LKSect2}
 L_K(u)(x)=\int_{\R^N}\delta (u,x,y) K(y)\d y, 
  \end{align}
  where 
\begin{equation*}
  	 \delta (u,x,y)=u(x+y)+u(x-y)-2u(x).
\end{equation*}
Equivalently, by the symmetry of the kernel,
  \begin{align}\label{LKpv}
  	L_K (u)(x)= 2 \,P.V.\ \int_{\rn} \{ u(y)-u(x)\} K(y-x) \d y ,
  \end{align}
  where PV stands for the principal value. 
Next, consider ${\cal M}^+$ and ${\cal M}^-$ from \eqref{Pucci}. They depend on the parameters $\Lambda$, $\lambda$ and $\alpha$, but we do not display this dependence for ease of notation. Explicitly,
  \begin{align}\label{def Pucci}
  	\mathcal{M}^+ (u)(x)=\int_{\rn} \frac{S_+(\delta(u,x,y))}{|y|^{N+2\alpha}} \d y , \;\; \mathcal{M}^- (u)(x)=\int_{\rn} \frac{S_-(\delta(u,x,y))}{|y|^{N+2\alpha}} \d y ,
  \end{align}
  where $S_+(t)=\Lambda t^+ - \lambda t^-$, and $S_-(t)=\lambda t^+ - \Lambda t^-$.
  \smallskip
  
We say $\mathcal F\in\mathcal E$ if $\mathcal F$ is either a Pucci extremal operator like $\mathcal F= \mathcal{M}^+$, $\mathcal F= \mathcal{M}^-$, or  an Isaacs type  integral operator in the form:
\begin{equation}\label{Isaacs}
	{\cal F}(u)=\inf_{a\in A}\sup_{b \in B}L_{K^{a,b}}(u),
\end{equation}
where $A$, $B$ are index sets, and for each ${a\in A}$ and ${b \in B}$, $L_{K^{a,b}}$ is a linear nonlocal operator for a kernel $K^{a,b}$ which is  $-(N+2\alpha)$ homogeneous.
Of course the study of Isaacs operators
  \begin{equation*}
  	{\cal F}(u)=\sup_{a\in A}\inf_{b \in B}L_{K^{a,b}}(u)
  \end{equation*}
 is completely analogous.

It is worth mentioning that our operators $\mathcal{F}\in \mathcal{E}$ are uniformly elliptic, i.e.\
  \begin{align}\label{SC}
  	\mathcal{M}^- (u-v)\leq \mathcal{F} (u)-\mathcal{F} (v)\leq \mathcal{M}^+ (u-v) ,
  \end{align}
  for all admissible functions $u$ and $v$. They are also scale invariant of order $2 \alpha$ in the sense that
  $\mathcal{F}(u_r)(x) =|r|^{2 \alpha}\mathcal{F} (u)(r x)$, where $u_r(x)=u(rx)$,
  for any $r \in \mathbb{R}$.\smallskip
  
  In what follows we recall some basic definitions and comparison theorems for $\mathcal{F}\in \mathcal{E}$ .
  \begin{definition}\label{DefVisc}
  	Assume $f:\mathbb{R}\rightarrow \mathbb{R} $ is continuous, and $g$ is a real function defined in $\rn$.
  	We say that $u$ is a viscosity super(sub)solution of
  	\begin{align}\label{def visc}
  		\mathcal{F} (u)+ f(u)=g(x) 
  	\end{align} 
  	at the point $x_0\in \rn$ where $u$ is continuous, if for any neighborhood $V$ of $x_0$ and for any $\varphi\in C^2(\overline{V})$ such that $u(x_0)=\varphi(x_0)$ and $u>\varphi$ $($resp. $u<\varphi )$ in $V\setminus\{x_0\}$, then 
  	\begin{align*}
  		\mathcal{F}(v)(x_0)+f(v(x_0)) \leq (\geq ) \; g(x_0) \; ,
  	\end{align*}
  	where $v:=u$ in $\rn \setminus V$, and $v:=\varphi$ in $V$.
  	We say that $u\in C(\Omega)$ is a viscosity super(sub)solution of \eqref{def visc} in $\Omega\subset \rn$ if it is a super(sub)solution of \eqref{def visc} at every point of $\Omega$.
  \end{definition}
  
  Now we recall the comparison principle for our operators in bounded domains. 
  \begin{proposition}\label{CP omega}
  	Assume $u$, $v\in C(\overline{\Omega})$ are bounded sub and supersolution satisfying 
  	$$
  	\mathcal F(u)\, \geq \, g \, \geq \, \mathcal{F} (v)\ \ {\rm in}\ \ \Omega,
  	$$
  	where $\Omega$ is an open bounded subset of $\mathbb R^N$ and $g$ is a continuous function in $\Omega$. If $u\leq v$ in $\mathbb R^N\setminus \Omega$, then $u\leq v$ in $\Omega$.
  \end{proposition}
  
  For the proof of Proposition \ref{CP omega} see \cite[Theorem 5.2]{CS1}. Notice that, when $v=0$, Proposition~\ref{CP omega} recovers the maximum principle for integral operators.
 
 Next we recall a strong maximum principle (SMP) in general cones stated in \cite[Theorem 2.8]{QSX}. 
 For reader's convenience we include a proof. We also need a boundary Harnack type result from \cite{RSpotential}, see also \cite{FRbook}.

  \begin{lemma}\label{SMP cone}
  \begin{enumerate}[(i)] Let $u$ be a solution of
  	\begin{equation*}
  		\left\{
  		\begin{array}{rlll}
  			{\mathcal F }(u)&\leq & 0\ \ \rm { in} \ \ \mathcal{C}_\omega , \\
  			u &\geq & 0\ \ \rm {in} \ \ \mathbb{R}^N.
  		\end{array}
  		\right.
  	\end{equation*}
  \item (Strong Maximum Principle) We have either $u>0$ in $\mathcal{C}_\omega$ or $u\equiv 0$ in $\mathcal{C}_\omega$.

\item (Half boundary Harnack)
Let  $\Omega \subset \mathbb{R}^n$ be any open set. Assume that there is $x_0 \in B_{1 / 2}$ and $\varrho>0$ such that $B_{2 \varrho}\left(x_0\right) \subset \Omega \cap B_{1 / 2}$.
Let $v\in C\left(B_1\right)$ be nonnegative such that 
$$
\left\{\begin{aligned}
	\M^{+}\left( v-au\right) &\ge 0  & \text { in } B_1 \cap \Omega \\
	v\ & = 0  \le u& \text { in } B_1 \backslash \Omega
\end{aligned}\right.
$$
for all $a\ge 0$.
Then,
$
v \le Cu$ in $ B_{1 / 2},
$
for a constant $C$ depending only on $N, \lambda,\Lambda,\alpha$, $\varrho$, and on the positive integrals 
$
\int_{\mathbb{R}^n} \frac{u(x)}{1+|x|^{n+2 s}} \d x
$
  and 
  $
  \int_{\mathbb{R}^n} \frac{v(x)}{1+|x|^{n+2 s}} \d x .
  $
  \end{enumerate}
  \end{lemma}
  
\begin{proof}
\textit{(i)}	Suppose there exists a $x_0\in \mathcal{C}_\omega$ such that $u(x_0)=0$. Then $L_{K}$ as in \eqref{LK} reads as
	$$
\textstyle 		L_{K}(u)(x_0)=\int_{\R^N}(u(x_0+y)+u(x_0-y))K(y)\d y.
	$$
	By \eqref{K} and $u\ge 0$  in $\mathbb{R}^N$, we get
	$$
\textstyle 	L_{K}(u)(x_0)\geq \lambda \int_{\R^N}(u(x_0+y)+u(x_0-y))\frac{\d y}{|y|^{N+2\alpha}}
\geq \lambda \int_{\mathcal{C}_\omega}(u(x_0+y)+u(x_0-y))\frac{\d y}{|y|^{N+2\alpha}}
	$$
for any even kernel satisfying \eqref{K}. Therefore, for $\mathcal{F}\in\mathcal{E}$ we get
	$$
\textstyle 		\mathcal{F}(u)(x_0)\geq \lambda \int_{\mathcal{C}_\omega}(u(x_0+y)+u(x_0-y))\frac{\d y}{|y|^{N+2\alpha}}.
	$$
	If $u(x_1)>0$ for some $x_1\in \mathcal{C}_\omega$, by continuity $u>0$ in $B_r(x_1)\subset {\mathcal{C}_\omega}$, then
	$$
\textstyle 		\lambda \int_{\mathcal{C}_\omega}(u(x_0+y)+u(x_0-y))\frac{\d y}{|y|^{N+2\alpha}}>c>0,
	$$
which contradicts the hypothesis ${\mathcal F }(u)\leq  0$ in $\mathcal{C}_\omega$. \smallskip
	
\textit{(ii)} This follows from the proof of \cite[Theorem 1.2]{RSpotential} by taking the functions $u_1=v /  \int_{\mathbb{R}^n} \frac{v(x)\d x}{1+|x|^{n+2 s}}  $ and $u_2=u / \int_{\mathbb{R}^n} \frac{u(x)\d x}{1+|x|^{n+2 s}} $.
\end{proof}

We finish the section with a theorem concerning solvability of Dirichlet nonlocal problems in bounded Lipschitz domains.
Let us consider the problem
\begin{align}\label{eq solvability}
	\textrm{$\mathcal{M}^\pm (u)  =0$ in $\Omega$, \;\; $u=g$\; in $\mathbb{R}^N\setminus \Omega$},
\end{align}
where $\Omega\subset \mathbb{R}^N$ is a bounded Lipschitz domain  and $g$ continuous in $\mathbb{R}^N\setminus \Omega$. As in \cite{AudritoRos}, we assume 
\begin{align}\label{phi1}
\textrm{$\quad|g (x)-g (z)| \leq C_0|x-z|^{\alpha_0} \quad$ for all $x \in \mathbb{R}^N \backslash \Omega, \; z \in \partial \Omega$},
\end{align}
for some constants $C_0$ and $\alpha_0 \in(0,1)$. In particular, $g$ is $C^{0, \alpha_0}$ on $\partial \Omega$, but not necessarily outside $\bar{\Omega}$. Besides, suppose the growth condition
\begin{align}\label{phi2}
|g(x)| \leq C_0\left(1+|x|^{\alpha_0}\right), \quad x \in \mathbb{R}^N \backslash \Omega .
\end{align}

\begin{proposition}\label{solvability nonlocal} Let $\Omega \subset \mathbb{R}^N$ be any bounded Lipschitz domain and $\alpha \in$ $(0,1)$.  
Then there exists a solution of \eqref{eq solvability}-\eqref{phi1}-\eqref{phi2}.
\end{proposition}
   
The above is a fully nonlinear counterpart of the viscosity solvability result in \cite[Corollary 4.2]{AudritoRos}. A complete proof for more general domains satisfying an exterior corkscrew condition can be found in details in the recent book \cite[Theorem 3.2.27]{FRbook}.

  \setcounter{equation}{0}
  \section{The Fundamental Solutions in Cones}\label{section fundamental}
  
  In this section we investigate fundamental solutions for extremal integral operators in cones. We first give the definitions of
  $\beta^-$ and $\beta^+$ respectively. We set
  $$
  \beta^-=\inf\{\beta<0\,|\,\textrm{ there exists }u\in H_\beta(\omega), \;u>0, \;{\cal F}(u)\leq 0\ {\rm
  	in}\ \mathcal{C}_\omega\}
  $$
  and
  $$
  \beta^+=\sup\{\beta>0\,|\,\textrm{ there exists }u\in H_\beta(\omega), \;u>0, \;{\cal F}(u)\leq 0\ {\rm
  	in}\ \mathcal{C}_\omega\},
  $$
  where 
  $$H_\beta(\omega)=\{u\in C(\overline{ \mathcal{C}_\omega}\setminus \{0\}) :u(tx)=t^{-\beta} u(x)\ {\rm for}\ x\in
  \mathcal{C}_\omega ,\ u(x)\geq 0\ {\rm for }\ x\in (\mathcal{C}_\omega)^c\},$$ 
  for any $\beta \neq 0$.
  Notice that a function in $H_\beta (\omega)$ is $-\beta$ homogeneous, and so determined by its values on $\omega$.
  Moreover, observe that
  \begin{align}\label{beta+(w0,w)}
  	\textrm{$\beta^+ (\mathcal{C}_{\omega_0})\geq \beta^+ (\mathcal{C}_\omega)$ \, and \, $\beta^- (\mathcal{C}_{\omega_0})\leq \beta^- (\mathcal{C}_\omega)$ \, for any $\emptyset\neq\omega_0\subset\omega$.}
  \end{align}
Indeed, since a positive supersolution in the cone $\mathcal{C}_\omega$ is also a positive supersolution in any smaller cone $\mathcal{C}_{\omega_0}$ it follows that $H_\beta (\omega_0)\supset H_\beta (\omega)$, so we use the sup and inf definitions of $\beta^+$ and $\beta^-$, respectively, to obtain \eqref{beta+(w0,w)}.
  
It is also important to have in mind that $\beta^-$ is increasing while $\beta^+$ is decreasing with respect to the operator, that is,
\begin{center}
 	$\beta^-(\mathcal{M}^-)\le \beta^-(\mathcal{F})\le \beta^-(\mathcal{M}^+)\le 0\le \beta^+(\mathcal{M}^+)\le \beta^+(\mathcal{F})\le \beta^+(\mathcal{M}^-) $.
\end{center}
  
In \cite[Lemmas 3.1 and 3.2]{ASS} it was obtained explicit examples which show that the set of test functions defining $\beta^\pm$ is nonempty and the corresponding set of admissible $\beta$ is bounded, via Pucci extremal local operators.
In particular, by using the approximation results in \cite{CS2} which describe nonlocal Pucci's operators through limits of local ones, one sees that the quantities $\beta^\pm$ are well defined when $\alpha$ is close to $1$. Other cases are analyzed in the next result.
  
We consider the dimension-like numbers $\tilde N^{\pm} > 0$ from \cite{FQ}, as mentioned in the Introduction, which satisfy
	$\M^{+} (|x|^{2\alpha - \tilde N^{+}})=0$  and $\M^{-} (|x|^{2\alpha - \tilde N^{-}})=0$ in $\rn$. 
  \begin{lemma}\label{lemma well defined}
(i) If $\tilde N^+ > 2\alpha$, then $\beta^+$ is well defined and $ \beta^+>0$. Also, $\beta^+>0$ when  $\mathcal{F}=\mathcal{M}^-$. 

(ii) $\beta^-$ is well defined and $\beta^-<0$ if either $\tilde N^-< 2\alpha$ or $\mathcal{C}_\omega \subset \{x\cdot \nu >0\}$, for some $\nu\in \mathbb{S}^N$.
  \end{lemma}
  
  \begin{remark}\label{remark H}
  	Note that  $\beta^+>0$ occurs when for instance $\mathcal{M}^-=-(-\Delta)^\alpha$ (when $\Lambda=\lambda$), for all $N\ge 2$. Also, $\beta^-<0$ is satisfied in particular when either $\alpha\in (0,1/2]$ for all $N\ge 1$, or when $\Lambda$ is close to $ \lambda$ (in this case we are close to the fractional Laplacian).
  \end{remark}
  
\begin{proof}
Observe that our construction of $H^\beta (\omega)$ makes every positive $-\beta$ homogeneus supersolution in $\mathbb{R}^N$ as a test function for the subsets defining $\beta^+$ and $\beta^-$. In particular, a fundamental solution in $\rn$ is a supersolution for the problem in $\mathcal{C}_\omega$. 
  	
(i) As far as the parameter $\tilde N^+ -2\alpha$ in Theorem 1.1 in \cite{FQ} is positive, then the definition of $\beta^+$ implies that $\beta^+ (\mathcal{F})\geq \tilde N^+ -2\alpha >0$ for any $\mathcal{F}\in \mathcal{ E}$. On the other hand, in the case of the operator $\mathcal{M}^-$ we have $\beta^+(\mathcal{M}^-)\ge \tilde N^- -2\alpha>0$ without requiring any hypothesis.
  	In particular, this is the case of the fractional Laplacian under $N\geq 2$, see Remark \ref{remark H}.
  	
(ii) Regarding $\beta^- <0$, it is enough to note that the function $u(x)= (x \cdot \nu)_+^\alpha$, for some $\nu\in \mathbb{S}^{N-1}$, solves $\mathcal{F} (u)\leq 0$ in $\{x \cdot \nu > 0\}$, which is positive and $\alpha$ homogeneous, see \cite{RSduke}.
\end{proof}

  Next, to fill in the range of constants $\alpha$, in what concerns the boundedness of $\beta^\pm$ in the general case, we need to work in an independent way, by exploiting the nature of the integral definitions near $-2\alpha$ and $N$. 
  This is the content of our next result.

  \medskip
  
  \begin{lemma}\label{claim beta+}
  	$\beta^->-2\alpha$ and $\beta^+<N$. 
  \end{lemma}
  
  \begin{proof}
Let $u\in H_\beta (\omega)$ be a positive function. For $\mathcal{F}\in \mathcal{ E}$ (see Section 2), we will analyze the integrability of the expression defining $\mathcal{F}( u )$ at a fixed point $x\neq 0$, according to $\beta$, as in the proof of Lemma~3.1 in \cite{FQ}. 
We are led to look at the expression $\int_{\mathbb{R}^N} \delta_\beta (y)K(y) \d y$, where  
\begin{center}
	$\delta_\beta (y)=\frac{ u(\frac{x+y}{|x+y|}) }{|x+y|^\beta }+\frac{ u(\frac{x-y}{|x-y|}) }{|x-y|^\beta }-\frac{ 2u(\frac{x}{|x|}) }{|x|^\beta } $,
\end{center}
foreshadowing a Pucci operator as in  \eqref{Pucci} or an Isaacs like operator as in \eqref{Isaacs} with $K=K^{a,b}$.

We start by fixing a connected subset $\omega_0\subset\subset \omega$ in $\mathbb{S}^{N-1}$ such that $x\in \omega_0$. Notice that $u$ is bounded from above and from below by positive constants in $\overline{\omega}_0$; namely $m_0=\inf_{\omega_0} u$, $M_0=\sup_{\omega_0} u$.  Thus, in view of homogeneity,  $u$ is bounded at infinity when $\beta>0$, whereas it is bounded at zero for $\beta<0$. 

We first consider $\beta\in (0,N)$. In this range the function $y\mapsto \delta_\beta (y)K(y)$ has three singularities: $y=0$, $y=x$ and $y=-x$.
Notice that the integral around $y=0$ is well defined independently on $\beta$ since $u$ is $C^2$ in $B_\eta(x)$ for $\eta$ small enough for $x\neq 0$ (or punctually $C^{1,1}$ at $x$ in the sense that there exists $M>0$ so that
$
|u(x+y)-u(x)-Du(x)\cdot y| \leqslant M|y|^2
$
for $y$ small). In fact, such regularity and the ellipticity condition \eqref{K} allows us estimate:
\begin{center}
	$
	\int_{B_\eta (0)} | \delta_\beta (y)|\, K(y) \d y \leq  2M  \int_{B_\eta (0)}  \frac{\d y}{|y|^{N-2+2\alpha}}<\infty,
	$
\end{center}
since $N-2+2\alpha<N$. We also say that the singularity at $y=0$ is removable in this case.

We now look upon the singularities $y=\pm x$. On the one hand,
\begin{center}
	$
	\frac{m_0 \lambda}{|x|^\beta } \int_{\mathbb{R}^N\setminus B_\eta(0)} \frac{ \d y }{|y|^{N+2\alpha} }\le
	\int_{\mathbb{R}^N\setminus B_\eta(0)} \frac{ u(\frac{\pm x}{|x|}) }{|x|^\beta }K(y)\d y \le 
	\frac{M_0 \Lambda}{|x|^\beta } \int_{\mathbb{R}^N\setminus B_\eta(0)} \frac{ \d y }{|y|^{N+2\alpha} }.
	$
\end{center}

On the other hand, regarding $y=x$, we may assume that  $0$ does not belong to the ball $B_{2\eta} (x)$ by making $\eta$ smaller if necessary. In particular, $B_\eta (x)\setminus B_\varepsilon(x)\subset \mathbb{R}^N\setminus B_\eta(0)$ and $|y|$ is comparable with $|x|$. More precisely, from $|y-x|\le \eta \le \frac{|x|}{2}$ one deduces
\begin{center}
		$\frac{|x|}{2} \le |x|-|x-y|\le |y|\le |y-x|+|x| \le \frac{3}{2}|x|$.
\end{center}
It is then easy to conclude that
\begin{align}\label{conta p2719}
\textstyle C_1\, \frac{\eta^{N-\beta}-\varepsilon^{N-\beta}}{N-\beta}\le
\int_{B_\eta (x)\setminus B_\varepsilon(x)} \frac{ u(\frac{x-y}{|x-y|}) }{|x-y|^\beta }K(y) \d y
 \le C_2 \,\frac{\eta^{N-\beta}-\varepsilon^{N-\beta}}{N-\beta} ,
\end{align}
when $\varepsilon<\mu$. Indeed, the positivity of $u$ and ellipticity of the kernel yield
\begin{center}
	$
\lambda\int_{B_\eta (x)\setminus B_\varepsilon(x)} \frac{ m_0}{|x-y|^\beta }\frac{\d y}{|y|^{N+2\alpha}} \le 
\int_{B_\eta (x)\setminus B_\varepsilon(x)} \frac{ u(\frac{x-y}{|x-y|}) }{|x-y|^\beta }K(y) \d y
\le \Lambda\int_{B_\eta (x)\setminus B_\varepsilon(x)} \frac{ m_0}{|x-y|^\beta }\frac{\d y}{|y|^{N+2\alpha}} .
$
\end{center}
Thus, \eqref{conta p2719} follows for $C_1, C_2 $  depending on $\lambda, \Lambda, m_0, M_0,\eta$ and $|x|$. 

Next, by passing the limit when $\eta,\varepsilon\to 0$ in \eqref{conta p2719} we end up concluding that the integral term comprising $|x-y|$ with singularity at $y=x$  is well defined whenever $\beta<N$.
Analogously we treat the singularity $y=-x$ appearing on the term  involving $|x+y|$.
Since $u$ is bounded in the remaining terms appearing in $\int_{\mathbb{R}^N\setminus B_\eta(0)} \delta_\beta (y)K(y)\d y,$ the operator of $-\beta$ homogeneous functions is always well defined for $\beta\in (0,N)$.

With respect to the limit as $\beta\to N$ in \eqref{conta p2719}, for $e =x|x|^{-1}$ we deduce
$$
\lim_{\beta\rightarrow N}\mathcal{F}( u(e)|x|^{-\beta} )=+\infty.
$$
In particular, $\mathcal{F}$ does not admit any positive $-N$ homogeneous supersolution, and so $\beta^+<N$.

\smallskip
  	
Finally, when $\beta\in (-2\alpha,0)$, the only singularity we have to check is at infinity. We observe that $\frac {|\delta_\beta(y)|}{|y|^{N+2\alpha}} \leq  \frac{C}{|y|^{N+2\alpha+\beta}}$, which is integrable outside any ball $B_{R}(0)$ since $N+2\alpha+\beta>N$.
More precisely, if $x\neq 0$ is fixed and $R$ is taken large enough so that $|x|\le \frac{R}{2}$, then $|x\pm y|\ge |y| -|x| \ge  |y|-\frac{R}{2} \ge \frac{|y|}{2}$ whenever $|y|\ge R$, yielding
  \begin{align*}
  \textstyle 	\int_{\rn\setminus B_R(0)} \delta_\beta(y) K(y)\d y 
  &   \textstyle 	\geq {2^{\beta +1} m_0 \lambda}	\int_{\rn\setminus B_R(0)} \frac{\d y}{|y|^{\beta+N+2\alpha}} - 	\frac{2M_0 \Lambda}{|x|^\beta } \int_{\mathbb{R}^N\setminus B_R(0)} \frac{ \d y }{|y|^{N+2\alpha} }\\
  &  \textstyle  \ge \frac{C_1}{\beta+2\alpha} \frac{1}{R^{\beta+2\alpha}}-C_2
  \end{align*} 
  for some positive constants $C_1$ and $C_2$ depending also on $|x|$.
Consequently, it blows up in the limit as $\beta\to -2\alpha$, from which 
  	$$
  	\lim_{\beta\rightarrow -2\alpha}\mathcal{F}( u( e)|x|^{-\beta} )=+\infty.
  	$$
Whence $\mathcal{F}$ does not have  $2\alpha$ homogeneous positive supersolutions, and $\beta^->-2\alpha$.
  \end{proof}

  In the sequel we prove some technical lemmas inspired in \cite{ASS} . The first one is a comparison principle for integral operators defined in cones, which is possible in these unbounded domains due to homogeneity of the involved functions.

  \begin{lemma}\label{CP}
  	Let $f\in H_{\beta+2\alpha}(\omega)$ nonnegative and $u,v\in H_\beta(\omega)$ satisfying
  	$$\left\{
  	\begin{array}{rllll}
  		-{\cal F}(u) \le f \leq -{\cal F}(v)&{\rm in} & \; \mathcal{C}_\omega 
  		\\ v >0 \,\,\;\qquad\qquad & {\rm in}&\; \mathcal{C}_\omega\\
  		u= 0 \,\,\;\qquad\qquad & {\rm in}&\, (\mathcal{C}_\omega)^c
  	\end{array}
  	\right.
  	$$
  	then either $u\leq v$ in $\mathcal{C}_\omega$ or there exists $t>1$ such that
  	$u=tv$.
  \end{lemma}
  
  \begin{proof}
  	The proof is similar to the proof of Proposition 2.5 in \cite{ASS}, by taking an extra care due to the nonlocal property of the operator. 
  	
  	Let $\omega_0 \subset\subset \omega$ connected. Denote $A_{r,s}$ the interception of the annulus $B_t \setminus \overline{B}_s$ with the cone $\mathcal{C}_\omega$, and $A^\prime_{s,t}$ the corresponding interception with $ \mathcal{C}_{\omega_0}$.
  	Set $a:=\frac{1}{2}\inf_{A^\prime_{1,2}} v>0 $, $C_0$ be such that $u\leq C_0$ in $A_{\frac{1}{2},4}$, and let $k>1$ be large enough such that $u\leq a k$ in $A^\prime_{\frac{1}{2},4}$. 
  	
  	Next, define $z=kv-u$, which is a viscosity solution of
  	\begin{center}
  		$-\mathcal{M}^-(z) \geq (k-1)f \geq 0$\; in $ A_{\frac{1}{2},4}$.
  	\end{center} 
  	Moreover, by construction,
  	\begin{center}
  		$z\geq ak$\; in $A^\prime_{{1},2}\,$, \, $z\geq -C_0$\, in $ A_{\frac{1}{2},4}\,$, \, and \, $z\geq 0$\, in $(\mathcal{C}_\omega)^c$.
  	\end{center} 
  	
  	\begin{claim}\label{lemma 2.4}
  		$z\geq 0$ in $A_{1,2}$\, up to taking a larger $k$.
  	\end{claim}

  	Denote $\Omega:= A_{\frac{1}{2},4}\setminus \overline{A^\prime}_{1,2}$. Define the functions $g_\pm \in L^1(\rn,({1+|x|)^{-N-2\alpha}})$ by
  	\begin{center}
  		$g_+ = 0 $ in $\rn \setminus A_{\frac{1}{2},4}\,$, \; $g_+=1$ in $\overline{A^\prime}_{1,2}\,$,
  	\end{center}
  	and 
  	\begin{center}
  		$g_-=1$ on $\mathcal{C}_\omega\cap(\partial B_{\frac{1}{2}}\cup \partial B_4)\,$,\; $g_-=0$ in $\partial \mathcal{C}_\omega\cap (B_3 \setminus B_{\frac{3}{4}})$,\; $g_-=0$ in $\overline{A^\prime}_{1,2}\,$,
  	\end{center}
  	where $0\leq g_\pm \leq 1$ in $\mathcal{C}_\omega$ are continuous, with $g_- (x) \leq 1/ (2|x|)^\beta$ in $\mathcal{C}_\omega \setminus A_{\frac{1}{2},4}.$ We also extend $g_- $ radially in $(\overline{\mathcal{C}_\omega})^c \cap B_4$, that is, $g_-$ is constant on the part of the spheres in $B_4$ which are contained in the exterior of the closure of the cone.
  	
  	Then, let $\psi_\pm$ be solutions of the following Dirichlet nonlocal problems
  	\begin{center}
  		$\mathcal{M}^\pm (\psi_\pm)  = 0 $\, in $\Omega$, \; $\psi_\pm = g_\pm$\, in $\rn \setminus \Omega$,
  	\end{center}
  	given by Proposition \ref{solvability nonlocal}, since $\Omega$ satisfies the exterior and interior uniform cone  property (see also the characterization of Lipschitz domains in \cite[Theorem 5.1.39]{teseLip}).
  	Thus, elliptic estimates and boundary Harnack (Lemma \ref{SMP cone}(ii), since $g_\pm =0$ in $(B_2\setminus B_1)\cap \mathcal{C}_\omega^c$) applied in the flat boundary yield some $\varepsilon_0 >0$, independent of $\psi_\pm$, such that 
  	\begin{center}
  		$\psi_+>\varepsilon_0\, \psi_-$\, in $\Omega\cap (B_2\setminus \overline{B}_1)=A_{1,2}\setminus \overline{A^\prime}_{1,2}$. 
  	\end{center}
  	
  	Hence, if $\psi = ak \psi_+ - C_0 \psi_- $, we have $z\geq \psi $ in $\rn \setminus \Omega$. So $z\geq \psi$ in $\Omega$ by Theorem~\ref{CP omega}. In particular, $z >0$ in $A_{1,2}\setminus \overline{A^\prime}_{1,2}$ provided $k> C_0 (\varepsilon_0 a)^{-1}$.
  	This proves Claim~\ref{lemma 2.4}.

  	Then, $z\geq 0$ in $\mathcal{C}_\omega$ by homogeneity. Next, we define the quantity
  	\begin{center}
  		$t:=\inf \{\, k>1, \, u\leq kv $\, in $\mathcal{C}_\omega \}$,
  	\end{center}
  	which is finite.
  	If $t=1$ then the first conclusion in Lemma \ref{CP} follows. 
  	
  	If $t>1$, we set $w=tv-u \geq 0$. 
  	In order to obtain a contradiction, assume $w\not\equiv 0$, then $w>0$ in $\mathcal{C}_\omega$ by SMP (Lemma \ref{SMP cone}(i)). 
  	Set $\tilde{z}=w-\delta v$, where $\delta\in (0,t-1)$ is small enough such that $\tilde{z}\geq 0$ in $ \mathcal{C}_{\omega_0}$, $\tilde{z}\geq a$ in $A^\prime_{1,2}$, and $\tilde{z} \geq - \varepsilon_0 \tilde{a}$ in $A^\prime_{\frac{1}{2},4}$, 
  	where $\tilde{a}=\frac{1}{2}\inf_{A^\prime_{1,2}} {w}$. 
  	Then, by repeating the argument in Claim \ref{lemma 2.4}, with $\tilde{z}$ in place of $z$, we obtain $\tilde{z}\geq 0$ in $\mathcal{C}_\omega$, a contradiction with the definition of $t$ as an infimum. This shows that $u=tv$ in $\mathcal{C}_\omega$.
  \end{proof}

  \begin{lemma}\label{le 2.2}
  	Assume that $\gamma>0$, $\beta\in (0, \beta^+]$, and $\psi\in H_\beta (\omega)$ is a nonnegative function. Set $\delta_\beta (y)= |e_1+y|^{-\beta}+|e_1-y|^{-\beta}-2$ and
  	$g(\beta)=\max\{c(\beta),-\beta\}$, where
  	\begin{equation}\label{2.5}
  		c(\beta) = \int_{\rn} \frac{S_+( \delta_\beta (y) )}{|y|^{N+2\alpha}} \d y .
  	\end{equation} 
  	Then there exists a unique nonnegative solution $u\in H_\beta(\omega)$ of
  	\begin{equation}\label{2.3}
  		\left\{
  		\begin{array}{rllll}
  			-{\cal F}( u)+g(\beta)(u-\psi)|x|^{-2\alpha}&=&|x|^{-2\alpha}(\beta \psi-\gamma u)+\beta |x|^{-\beta-2\alpha} &{\rm in}&\; \mathcal{C}_\omega \\
  			u &=&0 & {\rm in}& (\mathcal{C}_\omega)^c .
  		\end{array}
  		\right.
  	\end{equation}
  	Moreover, for some universal $c^*>0$, the following estimate holds
  	\begin{equation}\label{2.6}\|u\|_{L^\infty(\omega)}\leq
  		\frac{\beta}{\gamma}\,\{1+c^*\|\psi\|_{L^\infty(\omega)} \}.
  	\end{equation}
  	
  	Analogously, if $\gamma<0$, $\beta\in [\beta^-,0)$, and $\psi\in H_\beta (\omega)$ is nonnegative, then there exists a unique nonnegative solution $u\in H_\beta(\omega)$ of
  	\begin{equation}\label{2.3 phi-}
  		\left\{
  		\begin{array}{rllll}
  			-{\cal F}( u)+g(\beta)(u-\psi)|x|^{-2\alpha}&=&|x|^{-2\alpha}(\gamma u-\beta \psi)-\beta |x|^{-\beta-2\alpha} &{\rm in}&\; \mathcal{C}_\omega \\
  			u &=&0 & {\rm in}& (\mathcal{C}_\omega)^c .
  		\end{array}
  		\right.
  	\end{equation}
  	for $g(\beta)=\max \{c(\beta),\beta\}$, which satisfies the estimate \eqref{2.6} for some $c^*>0$.
  \end{lemma}
  
  \begin{proof}
  	We first notice that $0$ is a subsolution for the problem (\ref{2.3}). Indeed, if we replace $u=0$ into equation \eqref{2.3}, then the left
  	hand side is equal to 
  	$
  	-g(\beta) \psi |x|^{-2\alpha},
  	$
  	while the right hand side is equal to
  	$
  	\beta |x|^{-2\alpha}\psi +\beta |x|^{-\beta-2\alpha}.
  	$
  	
  	Now let us look at the function $$w_k(x)= k|x|^{-\beta}.$$
  	First notice that 
  	$ \mathcal{M}^+(w_k)$ given by \eqref{Pucci} is a well defined integral operator outside of $0$, since $y=0$ is a removable singularity; $y=x$ can be considered as a limit for $y$ in $\delta_1 \leq |y-x|\leq \delta_2$, and analogously for $y=-x$. 
  	Moreover, by a change of variable, 
  	\begin{align}\label{pucci+ wk}
  		\mathcal{M}^+(w_k)= { k c(\beta)}{|x|^{-\beta-2\alpha}}.
  	\end{align}
  	
  	Then we claim that $w_k$ is a supersolution of \eqref{2.3} for $k$ large enough. 
  	Indeed, it is enough to observe that $w_k$ satisfies 
  	$$ 
  	-\mathcal{M}^+(w_k)+g(\beta)(w_k-\psi)|x|^{-2\alpha}\geq |x|^{-2\alpha}(\beta \psi-\gamma w_k)+\beta |x|^{-\beta-2\alpha} \textrm{ in } \mathcal{C}_\omega ,
  	$$
  	once we chose
  	$$
  	k=\frac{\beta\,[1+c^*\sup_{\omega}\psi (x)]}{ g(\beta)-c(\beta) +\gamma},\quad \textrm{where\; $c^*=\sup_{0<\beta\leq \beta^+}\frac{g(\beta)+\beta}{\beta}$.}
  	$$
  	
  	\smallskip
  	
  	Next we define $u(x)=\sup_{v\in \mathcal{A}} v (x)$, where 
  	\begin{align*}
  		\mathcal{A}= \{ v\in C(\mathbb R^N\setminus\{0\}):\,
  		v & \textrm{ is subsolution of \eqref{2.3}}, \,v\leq w_k \textrm{
  			in } \mathcal{C}_\omega ,\, v\leq 0 \textrm{ in } (\mathcal{C}_\omega)^c \}.
  	\end{align*}
  	By using the preceding comparison lemma, and Perron's method for integral operators (see \cite{BCI} for instance),
  	we obtain that $u$ is a solution of the problem \eqref{2.3}. Now we observe that for any $v\in \mathcal{A}$, the function 
  	 $\tilde{v}$ defined as $\tilde v (x) = r^\beta v(r x)$ also belongs to $\mathcal A$, for any $r>0$. Thus, 
  	 \begin{center}
  	 	$u(rx) = \sup_{v\in \mathcal{A}} v (rx)= r^{-\beta} \sup_{\tilde v\in \mathcal{A}} \tilde v (x) =r^{-\beta} u(x)$,
  	 \end{center}
 that is, $u$ is $-\beta$ homogeneous. In addition, $u\geq 0$ in $\mathcal{C}_\omega$ since $0\in \mathcal A$. 
The next step is to show that this  solution is unique. 
  
  	Let $u_1 , \,u_2\in H_\beta (\omega)$ be solutions of \eqref{2.3}. 
  	We thereby claim that $\widetilde{u}_2 := u_2 +w_\varepsilon$ is a strict supersolution of \eqref{2.3}, for any $\varepsilon >0$.
  	To see this, notice that formally we have 
  	\begin{center}
  		$ \mathcal{F} (\widetilde{u_ 2}) -\mathcal{F}(u_2) \leq \mathcal{M}^+ (w_\varepsilon)< \varepsilon \{g(\beta)+\gamma\}|x|^{-\beta-2\alpha}$\; in $\mathcal{C}_\omega$,
  	\end{center}
  	where the latter comes from \eqref{pucci+ wk}. Then it is just a question of applying the definition of $u_2$ as a viscosity supersolution.
  	
  	In order to obtain a contradiction, assume $\varepsilon:=\max_{\overline{\omega}} (u_1-u_2) >0$.
  	By homogeneity of $u_1$ and $u_2$, we have that the strict supersolution $\tilde{u}_2$ touches $u_1$ by above at some point. But this is a contradiction with the definition of $u_1$ as a viscosity subsolution. 
  	Similarly we see that $\min_{\overline{\omega}} (u_1-u_2)<0$ derives a contradiction. Finally, the estimate \eqref{2.6} is a direct consequence of the choice of $k$. 
  	
  	If $\gamma,\beta<0$, then we still obtain a solution via Perron's method, which lies between $0$ and $w_k$, for the same $w_k$ as above; these are a pair of sub and supersolutions respectively for equation \eqref{2.3 phi-}.
  	In this case we choose
  	$$
  	k=\frac{-\beta\,[1+c^*\sup_{\omega}\psi (x)]}{ g(\beta)-c(\beta) -\gamma},\quad \textrm{with\; $c^*=\sup_{\beta^-\leq \beta\leq 0}\frac{g(\beta)-\beta}{-\beta}$.}
  	$$
  	This completes the proof of Lemma \ref{le 2.2}.
  \end{proof}

  \bigskip
  
  Now we are in the position of proving Theorem \ref{tm 0.1}. The proof of it is based on the following Leray-Schauder theorem due to Rabinowitz \cite{Ra}.
  \begin{proposition}\label{LS} Let $X$ be a real Banach space, $\mathcal  K\subset X$ a convex cone, and $A:[0,\infty)\times \mathcal  K\to \mathcal  K$ a completely continuous operator such that $A(0,u)=0$ for every $u\in \mathcal  K$. Then there exists an unbounded connected set $S\subset [0,\infty)\times \mathcal  K$ with $(0,0)\in S$, such that $A(\alpha,u)=u$ for every $(\alpha,u)\in S$.
  \end{proposition}
  
  \begin{proof}[Proof of Theorem \ref{tm 0.1}]
  	Let $\mathcal K $ be the convex cone of nonnegative continuous functions defined in the closure of the set $\omega$ in  $ {S}^{N-1}$. Given $\beta\in (0,\beta^+]$ and $u\in \mathcal K$, we consider $\psi$ as the extension of $u$ in a $-\beta$-homogeneous way to the cone, that is, $\psi (x) =|x|^{-\beta} u ( \frac{x}{|x|} )$.  For such $\psi$ and $\gamma=\beta^+$, let $\frak U\in  H_\beta(\omega)$ be the unique nonnegative solution of problem \eqref{2.3} given by Lemma \ref{le 2.2}, that is,
  	  \begin{equation*}
  		-{\cal F}( \frak { U})(x)+g(\beta)|x|^{-2\alpha}({\frak U} -\psi )(x)= |x|^{-2\alpha}(\beta \psi -\beta^+\,{\frak U} )(x)+\beta  |x|^{-\beta-2\alpha} \quad {\rm for}	\;\; x\in \mathcal{C}_\omega .
  	\end{equation*}
By the $2\alpha$ scale invariance of $\mathcal{F}$ and homogeneity of $\frak U$, the latter reads as
  \begin{equation*}
  	\textstyle	-{\cal F}( \frak { U})(e)=-g(\beta)({\frak U} - u) (e)+(\beta  u-\beta^+\,{\frak U} )(e)+\beta , \quad e=\frac{x}{|x|}.
  \end{equation*}
 Now, let $A(\beta,u)\in \mathcal  K$ be the restriction of ${\frak U}$ to $\overline{\omega}$. 
Also, set $A(\beta,u)=0$ for $u\in \mathcal  K$ and $\beta\leq 0$. 
  	Thus, Holder regularity estimates up to the boundary for integral operators in \cite{RSduke}, together with the estimate \eqref{2.6}, yield
\begin{center}
	  	$\| A(\beta,u) \|_{C^{0,\theta}(\overline \omega)}\le C (\beta, \beta^+, \|u\|_{L^\infty (\omega)}, N,\lambda,\Lambda) $,
\end{center}
  	for some $\theta \in (0,1)$.
Finally, nonlocal stability of viscosity solutions from \cite{CS2} implies that $A:[0,+\infty)\times \mathcal  K \rightarrow \mathcal  K$ is a completely continuous operator.
  	
  \smallskip	
  	Next, by the Leray-Schauder theorem (Proposition \ref{LS}), there exists an unbounded connected subset $S\subset[0,+\infty)\times \mathcal  K$ such that $(0,0)\in S$, and for any $(\beta, u)\in S$, we have $A(\beta,u)=u$.
  	In other words, for each $(\beta,\tilde u)\in S$, the function
  	$u(x)=|x|^{-\beta}\,\tilde u( {|x|^{-1}}x)$ is a solution of
  	\begin{equation}\label{2.4}
  		\left\{
  		\begin{array}{rllll}
  			-{\cal F}(u)&=&|x|^{-2\alpha}(\beta -\beta^+ )u+\beta |x|^{-\beta-2\alpha} &{\rm in}& \; \mathcal{C}_\omega \\
  			u&=&0 & {\rm in}& (\mathcal{C}_\omega)^c .
  		\end{array}
  		\right.
  	\end{equation}
  	
  	We infer that $S\subset [0,\beta^+]\times \mathcal  K$. Indeed, if $(\beta,\tilde u)\in S$ is such that $\beta>\beta^+$, then $u\in H_\beta (\omega)$ satisfies
  	\begin{align*}
  		-{\cal F}(u) \gneqq 0 \;\textrm{ in } \mathcal{C}_\omega \, , \;\;
  		u = 0 \; \; \textrm{ in }\; \mathbb \rn\setminus \mathcal{C}_\omega.
  	\end{align*}
  	Since $u\geq 0$ in $\mathcal{C}_\omega$ and the strict inequality holds above, we obtain $u>0$ in $\mathcal{C}_\omega$ by SMP; but this contradicts the definition of $\beta^+$. 
  	
  	Further, by the unboundedness of $S$, we deduce that for each $j\geq 1$, there exists $0\leq
  	\beta_j\leq \beta^+$ and $u_j\in H_{\beta_j}(\omega)$ such that
  	$\|u_j\|_{L^\infty(\omega)}\geq j$ and $A(\beta_j,\tilde u_j)=\tilde
  	u_j$, that is
  	$$
  	-{\cal F} (u_j)=|x|^{-2\alpha}(\beta_j -\beta^+ )u_j+\beta_j
  	|x|^{-\beta_j-2\alpha}.
  	$$
  	We assume that $\beta_j\to \tilde \beta\in [0,\beta^+]$, and we infer
  	that $\tilde\beta>0$. In fact, since $\|u_j\|_{L^\infty(\omega)}\rightarrow \infty$, we deduce from 
  	$$
  	j\leq \|u_j\|_{L^\infty(\omega)}\leq
  	\frac{\beta_j\,[1+c^*\|u_j\|_{L^\infty(\omega)}]}{\beta^+}
  	$$
  	that $\beta_j\geq c>0$. 
  	Next we define
  	$U_j=\frac{u_j}{\|u_j\|_{L^\infty(\omega)}}$, which is locally Holder continuous by \cite{CS2}. Then we may assume $U_j\to
  	U$ locally uniformly, for some nonnegative $\tilde{\beta}$ homogeneous function $U$ with $\|U\|_{L^\infty (\omega)}=1$. Moreover, $U$ satisfies
  	\begin{align}\label{eq MP}
  		-\mathcal F(U)=(\tilde\beta-\beta^+)|x|^{-2\alpha}U \;\textrm{ in } \mathcal{C}_\omega \, , \;\;
  		U = 0 \; \; \textrm{ in }\; \mathbb \rn\setminus \mathcal{C}_\omega .
  	\end{align}
  	Observe that $U>0$ in $\mathcal{C}_\omega$ by SMP, and so $U\in H_{\tilde\beta}(\omega)$. 
  	
  	Now we claim that $\tilde \beta = \beta^+$.
  	Otherwise, using $\tilde \beta<\beta^+$ in the definition of $\beta^+$, we find some  $\tilde \beta<\hat\beta\leq \beta^+$  and $v\in H_{\hat{\beta}} (\omega)$ such that $v>0$ and ${\mathcal{F}}(v)\leq 0$ in $ \mathcal{C}_\omega$. Now we define 
  	$$z=(v)^{\frac{\tilde \beta}{\hat \beta}},$$
 and by using the simple argument of Proposition 4 in \cite{Basic} we find 
  	that by concavity $$z(x)-z(y)\leq \frac{\tilde \beta}{\hat \beta} v(x)^{\frac{\tilde \beta}{\hat \beta}-1} (v(x)-v(y)) $$ 
  	this also holds for the test function and then we find that
  	$${\mathcal{F}}(z)\leq \frac{\tilde \beta}{\hat \beta} v(x)^{\frac{\tilde \beta}{\hat \beta}-1} {\mathcal{F}}(v)\leq 0\quad\mbox{ in }\quad C_\omega.$$
  	
  	Then applying Lemma~\ref{CP} to $U$ and to any positive multiple of $z$, we obtain that either $U=tv$ for some $t>0$, or $U\geq tv$ for every $t>0$. But the first alternative contradicts the strict inequality $\mathcal{F}(U)>0$ from \eqref{eq MP}; while the second leads to $U\leq 0$ in $\mathcal{C}_\omega$ which in turn contradicts $U>0$. So, $\tilde \beta = \beta^+$.
  	This also proves the uniqueness of $\beta^+$.
  	
  	Similarly one proves the existence of another fundamental solution which is $-\beta^-$ homogeneous, through the second part of Lemma \ref{le 2.2}. 
  \end{proof}

  \begin{remark}\label{def f}
  	Since the fundamental solutions $\phi^\pm$ are $-\beta^\pm$ homogeneous, they can be written as
  	$$
  	\phi^\pm(x)=\left\{
  	\begin{array}{ll}
  		{f^\pm(e)}{|x|^{-\beta^\pm}} &{\rm in}\quad \; \mathcal{C}_\omega \\ 0 & {\rm in}\quad (\mathcal{C}_\omega)^c ,
  	\end{array}
  	\right.
  	$$
  	where $e={|x|^{-1}}x$, and $f^\pm$ is the restriction of the fundamental solution $\phi^\pm$ to $\omega$.
  \end{remark}
  
  \section{The Liouville Theorem in Cones}\label{section Liouville}
  
  In this section, we prove Theorem \ref{tm 0.2} regarding a Liouville type theorem for extremal integral equation in cones.
  
  For fixed $\varepsilon>0 $, we define the function $w=w^\beta$ for any $\beta >0$ as:
  \begin{align}\label{def w}
  	w(x)=w^{\beta}(x) =\left\{
  	\begin{array}{lll}
  		{f(e)}{|x|^{-\beta}} &{\rm if}\quad |x|\geq\varepsilon, \quad x\in \mathcal{C}_\omega \\
  		{f(e)|x|}{\varepsilon^{-\beta-1}} &{\rm if}\quad |x|<\varepsilon, \quad x\in \mathcal{C}_\omega \\
  		0 &{\rm if}\quad \;\;x \in \rn \setminus \mathcal{C}_\omega\, ,
  	\end{array}
  	\right.
  \end{align}
  where $e={|x|^{-1}}x$, and for ease of notation $f=f^+=\phi^+|_{\,\omega}$ from Remark~\ref{def f}.
  
  Our next lemma settles $w$ as a subsolution in the cone $\mathcal{C}_\omega$ except for a ball centered at the vertex of this cone.
  
  \begin{lemma}\label{le 1.2}
  	Let $R_0$ be a fixed positive constant. Then for any $\beta\in (\beta^+ , N)$, there exists $\varepsilon_0=\varepsilon_0 (\beta,R_0)\in (0,1)$ such that, for each fixed $\omega_0\subset\subset\omega$ there exists $\rho=\rho (\omega_0,\varepsilon_0)$ such that
  	\begin{align}\label{F(w)>delta}
  		\mathcal{F} (w)(x)\geq \rho |x|^{-\beta-2\alpha} \textrm{ \;for all $x\in \mathcal{C}_{\omega_0}$ with $|x|\geq R_0$ }, \, \varepsilon\in (0, \varepsilon_0).
  	\end{align}
  \end{lemma}
  
  \begin{proof}
  	Fix some $\omega_0\subset\subset\omega$. 
  	We are going to show the more precise estimate \eqref{F(w)>delta} for some small $\rho=\rho (\omega_0,\varepsilon_0)$, by choosing $\varepsilon_0$ sufficiently small.
  	We start writing
  	\begin{align}\label{M- v}
  		\mathcal{F} \left(\frac{f(e)}{|x|^\beta}\right)=\frac{c(\beta)}{|x|^{\beta+2\alpha}},
  	\end{align}
  	for any $x\in \mathcal{C}_{\omega}$, $e={|x|^{-1}}x$, where $c(\beta)=c(e,\beta)$ is some function such that
  	$$ 
  	\intr \frac{S^-(\delta_\beta (y) )}{|y|^{N + 2 \alpha}}
  	\d y \leq c(\beta)\leq \intr \frac{S^+(\delta_\beta (y) )}{|y|^{N + 2 \alpha}}
  	\d y , \;\;
  	\delta_\beta (y)=\frac{f(\frac{e+y}{|e+y|})}{|e+y|^\beta}+\frac{f(\frac{e-y}{|e-y|})}{{|e-y|^\beta}}-2f(e).
  	$$
  	In fact, the regularity of $f$ and $\delta (v,x,0)=0$, for $v(x)=f(e)/|x|^\beta$, imply that $y=0$ is a removable singularity for the integral defined in \eqref{def Pucci}. The singularities $y=x$ and $y=-x$ are treated as limits when $y\in B_\epsilon \setminus B_\delta$ centered at $x$ and $-x$ respectively. So, \eqref{M- v} follows by a variable change, where $c(\beta)=\intr \frac{\kappa (\hat{y})}{|y|^{N + 2 \alpha}} \delta_\beta (y) \d y$, for some positive bounded function $\kappa\in [\lambda,\Lambda]$.
  	
  	We then infer that there exists $\beta_0> \beta^+$ such that
  	$c(\beta)\geq 2 \rho $ for all $\beta\geq \beta_0$ and $e\in \omega_0$. Indeed, a direct calculation shows that
  	$$
  	\frac{\partial^2 \delta_\beta}{\partial \beta^2}=
  	\frac{f(\frac{e+y}{|e+y|})\log^2|e+y|}{|e+y|^\beta}+\frac{f(\frac{e-y}{|e-y|})\log^2
  		|e-y|}{|e-y|^\beta} >0,
  	$$
  	and in particular the second order differential quotient of $\delta_\beta$ is positive in a small neighborhood of each point $\beta$. This allows us to use Fatou's lemma to pass the limit inside the integral, and conclude that $c(\beta)$ is strictly convex in $\beta$. 
  	Moreover, $c(\beta^-)=c( \beta^+)= 0$ for all $e\in\omega$. Therefore, since $c$ is increasing in $\beta\geq \beta^+$, we set $\rho= \frac{c(\beta_0)}{2}$ for some $\beta_0\in (\beta^+, N)$ in order to conclude the desired estimate on $c(\beta)$.
  	
  	Next we claim that, for $|x| \geq R_0$,
  	\begin{align}\label{tends to 0}
  	 \left|  \mathcal{M}^- \left( w-\frac{f(e) }{|x|^\beta} \right) \right| \le C \frac{\ep^{N-\beta}}{|x|^{N+2\alpha}}\,   .
  	\end{align}
  	To see this, we first choose $\ep\leq \frac {R_0}4$; then for $|x|\geq
  	R_0$, we have
  	\begin{eqnarray*}
  		&\,& \left| \mathcal{M}^- \left(w-\frac {f(e)}{|x|^\beta}\right) \right| \leq
  		\sup_{K}
  		\int_{B_\ep(-x)}\left| w(x+y)-\frac{f(\frac{x+y}{|x+y|})}{|x+y|^\beta}\right| K(
  		y)\d y\\&+&
  		\sup_{K} \int_{B_\ep(x)}\left| w(x-y)-\frac{f(\frac{x-y}{|x-y|})}{|x-y|^\beta}\right| K(y)\d y + \sup_{K} \int_{B_\ep(0)}\left| w(x)-\frac{f(e)}{|x|^\beta}\right| K(y)\d y ,
  	\end{eqnarray*} 
 for any positive even kernel $K$  satisfying \eqref{K}.
  	Now we observe that
  	$$
  	\int_{B_\ep(-x)}{w(x+y)}K( y)\d y\leq
  	C\frac{\ep^{N-\beta}}{|x|^{N+2\alpha}},
  	$$
  	$$
  	\int_{B_\ep(-x)}{\frac{f(\frac{x+y}{|x+y|})}{|x+y|^\beta}}K(
  	y)\d y\leq \frac C{|x|^{N+2\alpha}}\int_{B_\ep(0)}\frac{\d y}{|y|^\beta}=C\frac{\ep^{N-\beta}}{|x|^{N+2\alpha}}, \\ \vspace{0.1cm}
  	$$
  	and analogously for the other integrals. Since $\beta<N$, then $\ep^{N-\beta}\to 0$ as $\ep\to 0$, which proves the claim \eqref{tends to 0}.
  	Therefore, by taking a small $\varepsilon >0$, we obtain
  	\begin{align*}
  		\mathcal{F} (w)\geq \mathcal{F} \left(\frac{f(e)}{|x|^\beta}\right)+\mathcal{M}^- \left( w-\frac{f(e) }{|x|^\beta} \right) \geq \frac{c(\beta)}{|x|^{\beta+2\alpha}}-\frac{C}{|x|^{\beta+2\alpha}}\,  \ep^{N-\beta}|x|^{\beta-N} ,
  	\end{align*}
from where we conclude \eqref{F(w)>delta}.
  \end{proof}
  
  \medskip
  
  For the next lemmas it is convenient to introduce the following notation,
  $$
  \varphi_{s,t}^\beta(R)=\inf_{sR\leq |x|\leq tR,\, x\in
  	\mathcal{C}_\omega}\frac{u(x)}{w^\beta(x)},
  $$
  for $0<s<t<+\infty$, where $w^\beta$ is the function defined in \eqref{def w}.
  
  \begin{lemma}\label{lemma varphi>0}
  	Assume $\beta\in ( \beta^+,N)$, $R_1>\frac \ep s$, and let $u$ be a positive solution of $\mathcal{F}(u)\leq 0$ in $\mathcal{C}_\omega$. Then $\varliminf_{R\to \infty}\varphi_{s,t}^\beta (R)>0$. 
  \end{lemma}
  
  \begin{proof}
Let us fix $\beta$ and look at $\varphi^\beta_{s,t}$ for this $\beta$.
  	We are going to show that for any $\omega_0\subset\subset \omega$ we have the existence of a positive constant $C$ depending on $\lambda, \Lambda, N,\omega_0, \mathrm{dist}(\omega_0,\partial\omega)$ such that
  	\begin{center}
  		$u(x)\geq C w^\beta(x)$ for all $sR\leq |x|\leq t R$ with $x\in  \mathcal{C}_{\omega_0}$ and $R\geq R_1$.
  	\end{center}
  	
  	Denote $A_{r,s}(R)$ the interception of the annulus $B_{rR} \setminus \overline{B}_{sR}$ with the cone $\mathcal{C}_{\omega_0}$. \smallskip
  	
  	We first claim that $\varphi^\beta_{s,t}(R_1)>0$. 
Indeed, since we have a fixed radius $R_1$, we may apply Lemma \ref{SMP cone}(ii) twice to conclude that $u\ge c_1 \phi^- \ge c_2 \phi^+$ in a neighborhood of the flat boundary
 $\partial A_{s,r} (R_1)\cap \partial \mathcal{C}_{\omega}$, in which we also have $\phi^+(x)\ge  (R_1 +1)^{\beta-\beta^+} w^\beta$.
Now we use both the boundedness of $w^\beta$ and the positivity of $u$ in compact intervals of $\mathcal{C}_\omega$ to ensure that $u\ge c_3\, \omega^\beta$  in $A_{r,s}(R_1)$ for an appropriate $\omega_0\subset \subset \omega$.
Therefore we ensure the estimate $u\geq \tilde{c} w^\beta $ in the whole annular section $sR_1\leq |x|\leq tR_1$ of the original cone $ \mathcal{C}_\omega$, where $\tilde c = \min \{c_2 (R_1 +1)^{\beta-\beta^+} , c_3\}$.
In particular, $\varphi^\beta_{s,t}(R_1)\ge \tilde c>0$.

\smallskip
  		Consider the truncation $U$ of $u$ given by 
  	\begin{align}\label{truncation}
  		\textrm{$U=u$ in $\rn \setminus B_\varepsilon (0)\cup \{0\} $, \quad $U(x)=u(\varepsilon x/|x| )$ for $x\neq 0 $ in $B_\varepsilon(0)$.}
  	\end{align}
  	Further, set 
  	\begin{center}
  		$\displaystyle \phi_{s,t}^\beta (x)=
  		\inf_{sR\leq |x|\leq tR,\, x\in
  			\mathcal{C}_\omega}\, \frac{U(x)}{w^\beta(x)}$.
  	\end{center}
  	Notice that $\phi_{s,t}^\beta=\varphi_{s,t}^\beta$ for all $R\geq R_1$. 
In the sequel we are going to show that 
  	\begin{align}\label{eq F(U) small}
  		\mathcal{F}(U)\leq \rho|x|^{-\beta-2\alpha}\quad \mbox{ in } \rn \setminus B_R(0)
  	\end{align}
  	for large $R$.
  	To see this, we start splitting 
  	\begin{align*}
  		\int_{\R^N} &(U(y)-U(x)) K(y-x)\d y =\int_{\R^N\setminus B_\epsilon(0)}(u(y)-u(x)) K(y-x)\d y\\
  		&\quad+\int_{B_\epsilon(0)}(U(y)-u(x)) K(y-x)\d y \\
  		&=\int_{\R^N}(u(y)-u(x)) K(y-x)\d y+\int_{B_\epsilon(0)}(U(y)-u(y)) K(y-x)\d y,
  	\end{align*}
  	and so, up to a principal value sign from definition \eqref{LKpv},
  	\begin{align}\label{LK bound}
  		L_{K} (U-u) \le \int_{B_\epsilon(0)}(U(y)-u(y)) K(y-x)\d y.
  	\end{align}
Since $u\geq 0$, for any positive even kernel $K$ satisfying \eqref{K} we have
  	\begin{align}\label{eq critico rho}
  		\int_{B_\epsilon(0)}(U(y)-u(y)) K(y-x)\d y\leq \sup_{\partial B_\epsilon(0)}u\;\int_{B_\epsilon(0)}\frac{\d y}{|y-x|^{N+2\alpha}}\leq \frac{C\varepsilon^N}{|x|^{N+2\alpha}},
  	\end{align}
  	which in turn is less or equal than ${\rho}{|x|^{-\beta-2\alpha}}$ if we choose $\varepsilon$ small enough, for $|x|\ge R$ for $R\ge 1$.
  	Thus, by taking the supremum over $K$ in \eqref{LK bound} we obtain
  	\begin{align}\label{eq critico 2}
  		\mathcal{F}(U)-\mathcal{F}(u)\leq \mathcal{M}^+(U-u)\leq \rho |x|^{-\beta-2\alpha},
  	\end{align}
  	and by using $\mathcal{F}(u)\leq 0$ we deduce \eqref{eq F(U) small}.
  	
  	Next, Lemma \ref{SMP cone}(ii) as above implies that there exists $0<c_0<1$, such
  	that
  	\begin{align}\label{U>w c2}
  		\textrm{$U(x)\geq c_0\phi_{s,t}^\beta(R_1)w(x)$, \; for $x\in \mathcal{C}_\omega$ with $\ep\leq |x|\leq sR_1$,}
  	\end{align}
  	and moreover for $x\in \mathcal{C}_\omega$ with $|x|\leq \ep$ due to the truncation of $U$.
  	So \eqref{U>w c2} is true for $x\in \mathcal{C}_\omega$ with $ |x|\leq tR_1$.

  	Now, since $w(x)\to 0$ as $|x|\to \infty$, then for any $\mu>0$, we
  	can choose $\tilde R>0$, such that $\varphi_{s,t}^\beta(R_1)w(x)\leq \mu$
  	for $|x|>\tilde R$. 
  	For $C:=c_0\varphi_{s,t}^\beta(R_1)$, it follows that
  	$$
  	U(x)+\mu\geq C w(x),
  	$$
  	for $x\in\{x\,|\,x\in \mathcal{C}_\omega, |x|\leq tR_1\}\cup \{x\,|\,x\in \mathcal{C}_\omega,
  	|x|\geq \tilde R\}$.
  	
  	\medskip
  	
  	Finally, by $\beta>\beta^+$ and estimate \eqref{F(w)>delta} in Lemma \ref{le 1.2}, \begin{center}
  		$\mathcal{F}(\mu+U(x))\leq \rho |x|^{-\beta-2\alpha}\leq \mathcal{F}(C w)$ \; in $\{x\,|\,x\in \mathcal{C}_\omega,tR_1<|x|<\tilde R\}$
  	\end{center} for some small $\rho>0$. Then it follows by the comparison principle in bounded domains (Proposition~\ref{CP omega}) that
  	$$
  	\mu+U(x)\geq C w(x) \quad {\rm for\ all}\;\;
  	|x|\geq sR_1.
  	$$
  	Since $\mu$ is arbitrary, with $U=u$ and $\phi_{s,t}^\beta=\varphi_{s,t}^\beta$ for large $x$, it yields 
  	$u(x)\geq C w(x)$ for all $ |x|\geq sR_1$, as desired.
  \end{proof}

Next we show that the limiting case $\beta^+$ is still true in the sense that solutions of \eqref{eq 0.002} are comparable with $w^{-\beta^+}$ at infinity. Notice that $w^{\beta^+}$ coincides with the fundamental solution $\phi^+$ for large values of $|x|$, so we will be analyzing $\varliminf_{|x|\rightarrow \infty,\, x\in \mathcal{C}_{\omega_0}} \frac{u(x)}{\phi^+ (x)}$ for $\omega_0\subset\subset \omega$.

  \begin{proposition}\label{Hopf infty}
  	Let $\beta^+>0$ and $u$ be a positive solution of \eqref{eq 0.002} for $0<p\leq \frac{\beta^++2\alpha}{\beta^+}$. Then, any $\omega_0\subset\subset \omega$ there exists a positive constant $C$ such that
  	\begin{align*}
  		u(Re)\geq CR^{-\beta^+} \textrm{ \, for large $R$ and for all $e\in \omega_0$},
  	\end{align*}
  	where $C$ depends only on $\lambda,\Lambda,N$, $\omega_0$, and $\mathrm{dist}(\omega_0,\partial\omega)$.
  \end{proposition}
  
  \begin{proof}
Let $w=w^{\beta^+}$ as in (\ref{def w}) with $\beta=\beta^+$.
  	Since
  	$\mathcal{F}( {f(e) }{|x|^{-\beta^+}
  	} )=0$, it yields
  	\begin{align*}
  \mathcal{F}(w)(x)\ge \M^-(w-\phi^+) &= 
  \inf_K
  		 \int_{B_\varepsilon(x)\cup B_\varepsilon (-x)} \{ \delta(x,y,w)-\delta(x,y,f(e)|x|^{-\beta^+}) \}  K(y)\d y
  		 \\
  		 & \ge \inf_K
  		 \int_{B_\varepsilon(x)} A(y) K(y)\d y
  		 + \inf_K
  		 \int_{ B_\varepsilon (-x)} A(y)K(y)\d y, 
  	\end{align*}
  where
\begin{center}
$	A(y):=  \textstyle 
	w(x+y)+w(x-y)-f(\frac{x+y}{|x+y|})|x+y|^{-\beta^+}-f(\frac{x-y}{|x-y|})|x-y|^{-\beta^+}$.
\end{center}
\smallskip

Let $|x|>sR$ with $R>\frac{2\varepsilon}{s}$ and consider the integral over  $B_\varepsilon (x)$.
By the choice of $R$ we have $|x|>2\varepsilon$ so
 $|x+y|\geq |2x|-|y-x|\geq 3|x|/2\geq 3\varepsilon$, thus $w(x+y)=f(\frac{x+y}{|x+y|})|x+y|^{-\beta^+}$.
Now, since $K(y)\leq \Lambda |y|^{N+2\alpha}$, $f(\frac{x-y}{|x-y|})\leq M$, and 
 $\frac{|x-y|}{\varepsilon^{1+\beta^+}} \le 
 \frac{1}{|x-y|^{\beta^+}} $ in  $B_\varepsilon (x)$, we deduce
  	\begin{align*}
	 \int_{B_\varepsilon(x)} A(y) K(y)\d y & \ge
M\Lambda \int_{B_\varepsilon(x)}  \textstyle  \frac{1}{|y|^{N+2\alpha}}	 \left\{        \frac{|x-y|}{\varepsilon^{1+\beta^+}} - 
	 \frac{1}{|x-y|^{\beta^+}}
	 \right\}  \d y\\ \smallskip
	 &\ge   	-{ \textstyle \frac{2^{N+2\alpha}M\Lambda }{|x|^{N+2\alpha}}}
  		 \int_{B_\varepsilon (x)} \textstyle \frac{\d y}{|x-y|^{-\beta^+}}
  		=-C_1\frac{\varepsilon^{N-\beta^+}}{|x|^{N+2\alpha}},
  	\end{align*}
by using $|y|\ge |x| - |x-y|\ge |x| -\varepsilon \ge \frac{|x|}{2}$. A bound from below for the integral over $B_\varepsilon (-x)$ can be obtained in a similar way.
  	\smallskip
  	
Next, let us take $\omega_0\subset \subset \omega$. 
  	Since $0<p\leq \frac{\beta^++2\alpha}{\beta^+}$, we may choose some $\beta\in ( \beta^+,N)$ such that $\beta p<N+2\alpha$ and apply Lemma \ref{lemma varphi>0} to obtain
  	$$
  	\mathcal{F}(U)\leq -u^p(x)+\rho|x|^{-N-2\alpha}\leq - C_0\,|x|^{-p\beta}+\rho|x|^{-N-2\alpha}\leq -C_2|x|^{-N-2\alpha},
  	$$
  	in $sR\leq |x|\leq tR$, $x\in  \mathcal{C}_{\omega_0}$, for some small $\rho>0$.  Here, $U$ is as in \eqref{truncation} in Lemma \ref{lemma varphi>0} and the refined bound on $\mathcal{F}(U)$ is produced via \eqref{eq critico rho} into the estimate \eqref{eq critico 2}, where $\rho$ is chosen small enough.  By making $\varepsilon$ smaller, we may assume  $C_1\varepsilon^{N-\beta^+}\le {C_2}$, thus
  	\begin{center}
  		$\mathcal{F}(U)\leq \mathcal{F}(w)$ \;\; in $sR\leq |x|\leq tR$, $x\in  \mathcal{C}_{\omega_0}$. 
  	\end{center}
  	Using again that $w\to 0$ as $|x|\to +\infty$, the definition of the truncation $U$ and Hopf lemma, as in Lemma \ref{lemma varphi>0}, we have
  	\begin{center}
  		$U+\mu\geq w$ \;\; for \, $|x|\le sR$ \; or $|x|> tR$ when $x\in  \mathcal{C}_{\omega_0}$, \;or for $x\not\in  \mathcal{C}_{\omega_0}$,
  	\end{center}
  	for $\mu>0$.
  	Then, by Proposition \ref{CP omega} we deduce
  	\begin{center}
  		$U+\mu\geq w$ \;\; for \, $sR\leq |x|\leq tR$, \, $x\in  \mathcal{C}_{\omega_0}$. 
  	\end{center}
  	By letting $t\rightarrow \infty$ and $\mu\to 0$, we get
  	$u(x)\geq C R^{-\beta^+} $ for $ |x|\geq sR$,  $x\in  \mathcal{C}_{\omega_0}$,
  	which implies the desired result.
  \end{proof}
     
In the sequel we verify Theorem \ref{Hopf 0}. We carry over to a reflected version of our results by relying on the pure nature of the fundamental solution $\phi^-$ properly truncated at infinity. We skip the details and focus on the arguments that require a different analysis with respect to the case of $\phi^+$. 
  
  \begin{proof}[Proof of Theorem \ref{Hopf 0}]
  	For $R>0 $ and  $\beta<0$, let us denote the function $\widetilde{w}=\widetilde{w}^\beta$ as 
  	\begin{align}\label{def W}
  		\widetilde{w}(x)=	\widetilde{w}^\beta(x)=\left\{
  		\begin{array}{lll}
  			{g(e)}{|x|^{-\beta}} &{\rm if}\quad |x|\leq R, \quad x\in \mathcal{C}_\omega \\
  			{g(e)|x|^{-1} R^{1-\beta}} &{\rm if}\quad |x|> R, \quad x\in \mathcal{C}_\omega \\
  			0 &{\rm if}\quad \;\;x \in \rn \setminus \mathcal{C}_\omega\, ,
  		\end{array}
  		\right.
  	\end{align}
  	where $e={|x|^{-1}}x$, and  $g=f^-=\phi^-|_{\,\omega}$ from Remark~\ref{def f}. 
  For this proof we denote
  	\begin{equation}\label{phi- def}
  		\varphi_{s,t,\beta}^-\,(r): =\inf_{sr\leq |x|\leq tr, \, x\in \mathcal{C}_\omega} \,\frac{u(x)}{\widetilde{w}^\beta (x)} .
  	\end{equation}
  	
Let us fix a positive constant $\sigma$.
We first infer that a variant of Lemma \ref{le 1.2} shows that, for any $\beta\in (-2\alpha , \beta^-)$, there exists $R>2\sigma$ such that 
  	\begin{center}
  		${\cal F}(\widetilde{w}(x))\geq C_0 R^{-\beta-2\alpha}$\;\; for all $x\in \mathcal{C}_{\omega_0}$ with $|x|\leq \sigma$,\; $\sigma\leq \sigma_0$. 
  	\end{center}
  	Here $C_0$ is a constant depending on $\omega_0 \subset\subset \omega$, $\sigma_0$, $\alpha$, and $\beta$.
  	In turn, $C$ will be a constant which may change each step, with the same dependence as before.
  	
  	\smallskip
  	
  	Indeed, an analogous convexity argument proves that
  	\begin{center}
  		$\mathcal{F} ({g(e)}{|x|^{-\beta}} )={c(\beta)}{|x|^{-\beta-2\alpha}}\geq C_0 R^{-\beta-2\alpha}$ in $\mathcal{C}_{\omega_0}$\; for $|x|\leq \sigma$,
  	\end{center}
  	where $C_0 \geq C/\sigma_0$, by using that $|x|\leq \sigma\leq \sigma_0< \sigma_0 R/2$. 
  	W.l.g. we may assume $\sup_{\partial B_R (0)} u \, /(2\alpha +\beta)\leq C_0$ up to making $\sigma_0$ smaller such that $\sigma\leq \sigma_0$.\smallskip
  	
  	Next,
  	\begin{center}
  		$\mathcal{M}^- \left( \widetilde{w}-{g(e) }{|x|^{-\beta}} \right) \geq -C_0/2 \,{R^{-\beta-2\alpha}}$
  		for $|x|\leq \sigma$ in $\mathcal{C}_{\omega_0}$,
  	\end{center}
  	since the modulus of the LHS above is less or equal than
  	\begin{align*}
  		\sup_{K} \int_{\mathbb{R}^N} | \widetilde{w}(y)-g(y|y|^{-1}) |y|^{-\beta} | K(x-y) \d y 
  		&\leq \, 2 \sup_\omega g\, \sup_{K} \int_{\mathbb{R}^N\setminus B_R(0)} |y|^{-\beta} K(x-y) \d y \\
  		& \leq C \int_{\mathbb{R}^N\setminus B_R(0)} |y|^{-N-2\alpha-\beta}= C R^{-\beta-2\alpha},
  	\end{align*}
  	by using that $|x-y|\geq |y|-|x|\geq |y|/2$ and $\beta> -2\alpha$.\smallskip
  	
  	Moreover, by defining the truncation at infinity of $u$ given by $U=u$ in $B_R(0)$, and $U(x)= u(Rx|x|^{-1})$ if $|x|> R$, it yields
  	\begin{align*}
  		\mathcal{F}(U)-\mathcal{F}(u) \leq \mathcal{M}^+(U-u) \leq \sup_{K} \int_{\mathbb{R}^N\setminus B_R(0)} \{U(y)-u(y)\} K(x-y) \d y \leq C_0 R^{-\beta - 2 \alpha},
  	\end{align*}
  	by employing $u\geq 0$ and a bound of the truncation $U$ of $u$ on $\partial B_R (0)$. 
  	Thus, as in the proof of Lemma \ref{lemma varphi>0}, we obtain
  	$\varphi_{s,t,\beta}^-\,(r)  >0$ for $|x|\le \sigma$ and for $\beta \in (-2\alpha,\beta^-)$.
  \end{proof}

For the next lemma, let us consider the quotient
$$
\Phi_{s,t}^-(r):=\inf_{sr\le |x|\le tr, \, x\in \mathcal{C}_\omega}\frac{u(x)}{\phi^-(x)}.
$$

\begin{lemma} \label{lemma rho bounded}
	Let $\beta^-<0$ and $u$ be a positive solution of $\mathcal{F}(u)\le 0$ in $\mathcal C_\omega $ with $u\ge 0$ in $\mathbb{R}^N\setminus \mathcal C_\omega $. Then there exists a constant $C$ depending on $u$ and $r_0$ such that
	\begin{equation} \label{fundycmp-cones}
		0<\textstyle {\Phi}_{s,t}^-(r)\le C\quad \mbox{for large } \, r \ge \frac{t}{s} r_0.
	\end{equation}
	\end{lemma}

\begin{proof} Let $R>0$ fixed, $\widetilde{w}^{\beta^-}$ be the truncation \eqref{def W} with parameter $R$, and
	$$
	\hat{\varphi}^-(r):=\inf_{sr\le |x|\le tr, \,x\in  \mathcal{C}_\omega }\,\frac{u(x)}{\max\{0,\Psi (x)\}}\,, \quad \Psi (x):=\widetilde{w}^{\beta^-} - \sup_{sr_0\le |x|\le tr_0, \,x\in  \mathcal{C}_\omega} \widetilde{w}^{\beta^-}.
	$$
We claim that the map $ r \mapsto \hat {\varphi}^-(r)$ is nonincreasing in $ (\frac{t}{s} r_0,+\infty)$.
	Since $\beta^- < 0$, then 
	\begin{center}
		$\sup_{sr\le |x|\le tr, \,x\in  \mathcal{C}_\omega} \widetilde{w}^{\beta^-} > \sup_{sr_0\le |x|\le tr_0, \,x\in  \mathcal{C}_\omega} \widetilde{w}^{\beta^-}$ \; for $r>\frac{t}{s} r_0$.
	\end{center} Thus for $r> \frac{t}{s}r_0$ the quantity $\hat\varphi^-(r)$ is finite and positive, with
	$u +\mu \geq \hat\varphi^-(r) \widetilde{w}^{\beta^-}$ when either $|x|\le  tr_0$ or $|x|\ge tr$,  by construction of $\Psi$, and the definition of $\tilde{w}^{\beta^-}$ as a truncation for all $|x|> R$, which converges to $0$ as $|x|\to +\infty$.
	By the comparison principle we then have
	\begin{equation*}
		u(x)+\mu \geq \hat\varphi^-(r) \Psi (x) \quad \mbox{for} \;\;  t r_0 \le |x|\le tr, \;x\in \mathcal  \mathcal{C}_\omega
	\end{equation*}
	for $r\in ( \frac{t}{s}r_0, R)$, and the claim follows by letting $\mu \to 0$. Finally, $\Phi_{s,t}^-(r) \leq \hat\varphi^-(r)\le \hat\varphi^-(\frac{t}{s} r_0)$, for large $r$. Since $R$ is arbitrary the lemma is proven.
\end{proof}

Given $\omega_0\subset\omega$ such that $\overline{\omega}_0\subset \omega$ and
define, for $u_r(x):=u(rx)$ and for $r\ge 2 r_0$,
\begin{center}
	$	m(sr,tr,\omega_0):=\inf_{sr\leq |x|\leq tr,\, x\in  \mathcal{C}_{\omega_0}} u(x)=\inf_{s\leq |x|\leq t,\, x\in  \mathcal{C}_{\omega_0}} u_r(x)$.
\end{center}	
	
Notice that  Proposition \ref{Hopf infty} says that $	m(sR,tR,\omega_0)\geq {C}{R^{-\beta^+}}	$ for $R$ large enough.
In order to have estimates also from above over the infimum of $u$ in annular portions of smaller cones, we need the next lemma.

\begin{lemma}\label{lemaQ1}
Let $p\in \mathbb{R}$ and $u$ be a supersolution of \eqref{eq 0.002}.	For any $0<s_1<s_2<t_2<t_1<\infty$ and $\omega_2\subset \omega_1
\subset \omega$, we have
	$$
	\inf_{s_1R<|x|<t_1R,\, x\in  \mathcal{C}_{\omega_1} } u^p (x)
	\leq
	CR^{-2\alpha}\,
	m(s_2R, t_2R,\omega_2),
	$$
where $C$ depends on $\lambda,\Lambda, N, \omega_i, s_i, t_i$, $i=1,2$.
In addition, by setting $m_i(R):=m_i (s_i R, t_i R,\omega_i)$ for $i=1,2$, it follows:
\begin{enumerate}[(i)]
	\item if $p>0$ then $m_1(R)^p\le R^{-2\alpha} m_2(R)$ for large $R$.
	
	\item if $p\le 0$ then $m_2(R)\ge R^\frac{2\alpha}{1-p} $ for large $R$. In particular, any supersolution $u$ of \eqref{eq 0.002} is unbounded at infinity.
\end{enumerate}
\end{lemma}

\begin{proof}
	Let us choose a cut-off function $\rho$ with values between $0$ and $1$, such that
	$$ \rho(x)=\left\{
	\begin{array}{ll}
		1 &{\rm if}\quad x\in \mathcal{C}_{\omega_2}\ {\rm and}\ s_2<|x|<t_2, \vspace{0.1cm}
		\\ 0 &{\rm if}\quad |x|\geq t_1\ {\rm or}\ |x|\leq s_1\ {\rm or}\ x\not\in  \mathcal{C}_{\omega_1} \\
	\end{array}
	\right.
	$$
	and $\mathcal{M}^- (\rho(x))\geq -C$.
Now we define $\eta(x)=m_2(R)\rho(\frac x R)$ and
	$\xi(x)=u(x)-\eta(x)$. \smallskip
	
Set $\Omega_i :=\{x\in  \mathcal{C}_{\omega_i} :
	s_iR<|x|<t_iR\}$, $i=1,2$.
Obviously we have $\xi(x)=u(x)>0$ in
	$\partial\Omega_1 $. Also,  there exists $\bar x\in \bar \Omega_2$ such that $u(\bar
	x)=\inf_{x\in \Omega_2} u(x)$. So $\xi$ must possess a minimum $x^*\in
	\Omega_1 $. 
Note that $\xi (x^*)\le 0$, so in particular $u(x^*)\le m_2(R) $. 
	
	\smallskip
	
	It follows from the definition of viscosity solution and ellipticity \eqref{SC} that
	\begin{eqnarray*}
		-u^p(x^*)\geq \mathcal F (\eta (x^*))\geq \mathcal{M}^- (\eta(x^*))
		= \inf_{\Omega_2} u \;
		\mathcal{M}^- (\rho( {x} /R ))|_{x=x^*} ={R^{-2\alpha}}\inf_{\Omega_2} u \;\mathcal{M}^- ( \rho({x^*} )) .
	\end{eqnarray*}
	
	Hence, we conclude that
	$
\inf_{\Omega_1} u^p
\leq u^p(x^*)\leq C{R^{-2\alpha}}
\inf_{\Omega_2} u .
	$ \smallskip
	
We observe that for $p\le 0$ we have $m_2(R)^p\le \inf_{\Omega_1} u^p
\leq u^p(x^*)\leq C{R^{-2\alpha}}
m_2(R)$, which gives the statement. Meanwhile, for $p>0$ we deduce the conclusion \textit{(i)}.
\end{proof}

\begin{proposition}\label{propPhi2leC} Let $\beta^-<0$ and $u$ be a supersolution of \eqref{eq 0.002} with $p\ge  \frac{\beta^-+2\alpha}{\beta^-}$.
Given $0<s<t$ and $\omega_0\subset \subset \omega$, then
\begin{align}\label{est D}
	{\textrm{$\Phi_{s,t,\omega_0}^{-}(R):=\inf_{sR\le |x|\le tR, \, x\in \mathcal{C}_{\omega_0}}\frac{u(x)}{\phi^-(x)}\le C$ \; for large $R$}},
\end{align}
for a constant $C$ which depends on $t-s$ and $\omega\setminus \omega_0$.
\end{proposition}

The difference of Proposition \ref{propPhi2leC} with respect to Lemma \ref{lemma rho bounded} is that, since $\omega_0 \neq \omega$, then $\Phi_{s,t,\omega_0}^{-}$ could be larger than the function $\Phi_{s,t}^{-}$ there. As in \cite{AScpde}, a proper use of the comparison principle helps us to overcome this difficulty.

In addition, note that Proposition \ref{propPhi2leC} complements Proposition \ref{Hopf infty} in the sense that, whenever $p\in [p^-_*,p_*^+]$, where $p^\pm_*= \frac{\beta^\pm +2\alpha}{\beta^\pm }$, the infimum over compact sets of the cone of any nonnegative nontrivial supersolution of \eqref{eq 0.002} needs to be controlled at infinity by the fundamental solutions $\phi^+$ and $\phi^-$.

\begin{proof} Of course the conclusion of the lemma is trivial if $m(sR,tR,\omega_0)$ is bounded. So, we only have something to prove when $m(sR,tR,\omega_0)$ is an unbounded function of $R$.
	
Next, as far as $m(sR,tR,\omega_0)\to \infty$ then, up to taking a smaller cone and $s_0>s$, $t_0<t$ if necessary, we may assume that $m(sR,tR,\omega_0)\ge C_0R^{-\beta^-}$ for large $R$. Indeed, if $\omega_2\subset\subset \omega_1 \subset\subset\omega$, $s_1<s_2<t_2<t_1$, and $m_i$ are as in the proof of the Lemma \ref{lemaQ1}, we deduce
\begin{center}
	$C R^{-2\alpha}m_2(R)\ge u^{p-p^-_*}(x^*)  u^{p^-_*}(x^*)\ge m_1(R)^{p-p^-_*}m_2(R)^{p^-_*}$.
\end{center}
Thus, if $m_1(R)\to \infty$ as $R\to\infty$, in particular $m_1(R)\ge C$ for $R\ge R_1$, and we get
\begin{center}
	$m_2(R)\ge C R^{\frac{2\alpha}{1-p^-_*}}=CR^{-\beta^-}$ for $R\ge R_1$.
\end{center}
		
Set
\begin{equation*}
	w(x): = u(rx) - \delta\, \Phi_{s,t,\omega_0}^{-}(r) \,\phi^-(rx),
\end{equation*}
for $r> \frac{t}{s}r_0$, and $0 < \delta < 1$ to be chosen.  We assume $s=1$ and $t=2$ for the sake of simplicity and so omit the subindexes $s,t$ in the notation, namely $\Phi^-$ and $\Phi_0^{-}$. Then, for large $R$,
\begin{center}
	$-\mathcal{F} (w) \geq 0$ in $(B_4\setminus B_{\frac{1}{2}})\cap \mathcal{C}_\omega$, \; $ w \geq 0 $ on $( B_4\setminus B_{1/2})\cap \partial \mathcal{C}_\omega$ and in $\mathbb R^N\setminus \mathcal C_\omega$,\medskip
	
	$w \geq -C \delta  \,\Phi^{-}_0(r)\,  r^{-\beta^-}$ in $B_4\cap \mathcal C_\omega,$  \; $w \geq c (1-\delta)  \Phi^{-}_0(r)\,  r^{-\beta^-}\ge c_0$ in $(B_2\setminus B_1)\cap \mathcal C_{\omega_0}$.\smallskip
\end{center}

Set $\Omega : = (B_4\setminus B_{\frac{1}{2}} ) \cap \mathcal{C}_{\omega} \setminus ((B_2\setminus B_{1})\cap \mathcal{C}_{\omega_0})$ and let $v_1$ and $v_2$ be solutions of the following nonlocal Dirichlet problems
\begin{center}
	$-\mathcal{F}(v_i) = g_i$ \; in $\Omega$,\quad $v_i= g_i$ \; on $\partial\Omega$.
\end{center}
where $g_1 = g_2 = 0$ on  $ \partial\Omega\cap \partial \mathcal C_\omega $ and in $\mathbb{R}^N \setminus\mathcal{C}_\omega$. Also, $g_1= c (1-\delta)  \Phi^{-}_0(r)\,  r^{-\beta^-}>0$ and $g_2 = 0$ on the inner boundary of $\Omega$ as well as in its interior $(B_2\setminus B_1)\cap \mathcal C_{\omega_0}$. In addition, $g_1=0$ and $g_2 = C  \,\Phi^{-}_0(r)\,  r^{-\beta^-}$ on the top boundary part of $\Omega$ given by $\partial B_4 \cap \mathcal C_\omega $, as well as in $\overline B_{\frac{1}{2}}\cap \mathcal C_\omega$. In $(\mathbb R^N\setminus B_4)\cap \mathcal C_\omega$ we define $g_1$ as zero and $g_2$ extended in an integrable way.

Now, elliptic estimates and the boundary Harnack type result in Lemma \ref{SMP cone}(ii) applied on the flat boundary yield, for $\ep > 0$ sufficiently small, that
\begin{equation*}
	v_1 > \ep v_2 \quad \mbox{in } \ \Omega \cap \left( B_2 \setminus B_1\right) =  \left( B_2 \setminus B_1 \right)\cap \left( \mathcal{C}_\omega \setminus \mathcal C_{\omega_0} \right).
\end{equation*}
Set $v: = v_1 - \ep v_2$ and observe that
$w \geq v$ in $\mathbb R^N\setminus \Omega$ by choosing $\delta\le  \varepsilon$.
Hence, by the comparison principle we have $w \geq v$ in $\Omega$. In particular, $w \ge 0$ in $\left( B_2 \setminus B_1 \right)\cap \left( \mathcal{C}_\omega \setminus \mathcal C_{\omega_0} \right)$. 	
In other words, we have
\begin{equation*}
	u(rx) \geq \delta \,\Phi^{-}_0(r) \phi^-(rx)  \quad \mbox{for } \ 1\le |x|\le  2 , \, x\in \mathcal C_\omega.
\end{equation*}
Therefore
$
\Phi^{-}(r) \geq \delta \Phi^{-}_0( r).
$
Since $\Phi^-(r) \leq C$ by Lemma \ref{lemma rho bounded}, and $\phi^-$ is homogeneous, the proof is finished.
\end{proof}

\begin{proof}[Proof of Theorem \ref{tm 0.2}.]
Let $p>0$.
With respect to $m_i$ defined in Lemma \eqref{lemaQ1}, two situations might occur: either 
(i) $m_i(R)\le C_0 $ all subsets $\omega_i\subset \subset \Omega$; or (ii) there exists some $\omega_i$ such that  $m_i(R)\to +\infty$ as $R\to +\infty$ for some subset $\omega_i\subset \subset \omega$. 

We rule out possibility (ii) when $p>0$. Up to consider a smaller cone $\omega_1$, if $m_1(R)\to +\infty$ as $R\to\infty$, then by Lemma \ref{lemaQ1} and Proposition \ref{propPhi2leC} we get 
\begin{center}
	$C\le m_1(R)^p \le Cm_2(R)R^{-2\alpha}\le CR^{-(\beta^- +2\alpha)} \to 0$ as $R\to \infty$,
\end{center}
which is impossible.

Hence, (i) is in force.
Next we prove that for any $q\in(\frac{2\alpha} p,
\frac{2\alpha}{p-1})$ if $p>1$, while $q\in (\frac{2\alpha} p,+\infty)$ if $p\le 1$, given $0<s<t<\infty$, we have the more precise estimate
\begin{align}\label{eqleQ2}
	\textrm{	$m_1(R)\leq {C}{R^{-q}}$ \quad for large $R$.}
\end{align}

To see \eqref{eqleQ2}, we denote by $\gamma_0=\frac{2\alpha}{p}$ and $\gamma_i=\frac{2\alpha+\gamma_{i-1}}{p}$. It is easy to check that $\gamma_i\to \frac{2\alpha}{p-1}$ as $i\to \infty$ if $p>1$, and it goes to infinity if $p\leq 1$.
So there exists $k>0$ such that $\gamma_{k-1}>q$.
Then we pick up sequences
\begin{center}
	$
	s_1<s_2<...<s_{k-1}<s_k<t_k<t_{k-1}<...t_2<t_1
	$
	, \;	$
	\emptyset\neq \omega_k\subset \omega_{k-1}\subset...\subset \omega_1.
	$
\end{center}
By Lemma \ref{lemaQ1}(i) we know that
$
m_i^p(R)\leq {C_i \, m_{i+1}(R)}{R^{-2\alpha}},
$
where $m_i(R)= m(s_iR, t_iR,\omega_i)$ and $C_i$ is independent of $R$.
Since the function $m_i$ is bounded at infinity for any $i$, then we have
$$
m_{k-1}^p(R)\leq C_{k-1}\frac{m_k(R)}{R^{2\alpha}}\le \frac{M
	C_{k-1}}{R^{2\alpha}},
\;\;\;\textrm{ or } \;\;
m_{k-1}(R)\leq \frac{b_{k-1}}{R^{\frac {2\alpha} p}}.
$$
Similarly we deduce
$$
m_{k-2}(R)\leq \frac{b_{k-2}}{R^{\frac {2\alpha}p+\frac{2\alpha}{p^2}.
}}
$$
By iterating the above procedure, for $R$ large enough, we obtain
$
m_1 (R)\leq \frac{b_1}{R^{\gamma_{k-1}}}\leq \frac{b_1}{R^q}$, which gives us \eqref{eqleQ2}.

\smallskip

Let $p\in (0,\frac{\beta^++2\alpha}{\beta^+})$. Then, by combining \eqref{eqleQ2} and Proposition \ref{Hopf infty} we get 
\begin{center}
$cR^{-\beta^+}\le m_1(R)\le C R^{-q} $ \quad for large $R$
\end{center}
where $q\in (\frac{2\alpha}{p}, \frac{2\alpha}{p-1})$ if $p>1$ and $q\in (\frac{2\alpha}{p},+\infty)$ if $p\in (0,1]$.

\smallskip

So, if $p>1$, then we choose $q=q(p):=\frac{2\alpha}{p-1}-\epsilon_p$ with $\epsilon_p$ small enough such that $\epsilon_p< \frac{(p^+_*-p)\beta^+}{p-1}$. Recall that $p^+_*=\frac{\beta^++2\alpha}{\beta^+}$. Thus,
\begin{center}
$c\le CR^{-\left\{\frac{2\alpha}{p-1}-\beta^+ \right\}} R^{\epsilon_p}
=  C R^{-\frac{\beta^+}{p-1} (p^+_*-p)  }R^{\epsilon_p}\to 0 $ \quad  as $R\to\infty$
\end{center}
which is impossible. 

On the other hand, if $p\in (0,1]$ it is enough to take any $q=q(p) >\beta^+$ in order to obtain the absurdity
$c\le CR^{\beta^+-q} \to 0$ as $R\to\infty$. This concludes the proof in the positive subcritical case.

\smallskip

Now, let $p=\frac{\beta^++2\alpha}{\beta^+}$. We then prove the refined log-estimate
\begin{align}\label{eq inf log}
	u(x)\geq C\mathrm{log}(1+|x|)|x|^{-\beta^+}\; \textrm{ for } s_1R \le |x| \le r_1 R, \; x\in \mathcal{C}_{\omega_0} , \; R\ge R_1.
\end{align}

We split the proof of \eqref{eq inf log} in two cases. Firstly we assume $\beta^+\in (1,N)$ and set
\begin{align*}
 W(x)=\mathrm{log}(1+|x|)w(x), \quad |x|\geq s_1R_1> 2\varepsilon ,
\end{align*}
where $e=x/|x|$, $w=w_{\beta^+,\,\varepsilon}$\, from \eqref{def w} that is,
\begin{align*}
	w(x)=\left\{
	\begin{array}{lll}
		{f(e)}{|x|^{-\beta^+}} &{\rm if}\quad |x|\geq\varepsilon, \quad x\in \mathcal{C}_\omega \\
		{f(e)|x|}{\varepsilon^{-\beta^+-1}} &{\rm if}\quad |x|<\varepsilon, \quad x\in \mathcal{C}_\omega \\
		0 &{\rm if}\quad \;\;x \in \rn \setminus \mathcal{C}_\omega\, ,
	\end{array}
	\right.
\end{align*}
and $\phi^+(x)=\frac{f(e)}{|x|^{\beta^+}}$ satisfies $\mathcal{F}(\phi^+)=0$.
Thus for fixed $x\in \mathbb{R}^N\setminus B_{s_1R_1}$ we have
$$
\mathcal{F}(w)(x)\ge \inf_{b}\intr \delta(x,y,\phi^+)
\,K^b (y)
+\inf_b
\int_{B_\varepsilon(x)\cup B_\varepsilon(-x)} \{\delta (x,y,w)-\delta (x,y,\phi^+)\}\,K^b(y),
$$
since $|x|>\varepsilon$, and therefore 
$\mathcal{F}(w)\ge  -C|x|^{-\beta^+-2\alpha}.$
Now we infer that 
\begin{center}
	$\mathcal{F}(W)\geq - C|x|^{-\beta^+-2\alpha}$\; in $\mathbb{R}^N\setminus B_{s_1R_1}$.
\end{center}
Indeed, this follows by the proof of Lemma 6.1 in \cite{FQ}, by considering instead of
$\delta_1$ there,
\begin{center}
	$\delta_1(r,z)=g(|e+z|,\theta)f(\frac{e+z}{|e+z|})+g(|e-z|,\theta)f(\frac{e-z}{|e-z|})$,
\end{center} 
with $e_1 $ there replaced by $e$, since in our case we do not have a radially symmetric function. The analysis is carried out in the same way in $g$, and using the boundedness of $f$ in $\omega_0$.

We consider the truncation $U(x)$ defined in \eqref{truncation}. Then, as in \eqref{eq F(U) small} one finds that 
\begin{center}
	$\mathcal{F}(U)+U^p \le \rho |x|^{-\beta^+ -2\alpha}=\rho |x|^{-p\beta^+} $ \; for $|x|\ge s_1R$,
\end{center} 
up to taking larger $s_1$ if necessary.

Now, by Proposition \ref{Hopf infty} one finds $u(x)\ge C_0|x|^{-\beta^+}$ in $(B_{t_1R}\setminus B_{s_1R})\cap \mathcal{C}_{\omega_0}$, so
\begin{center}
	$-\mathcal{F}(W )\leq  C|x|^{-\beta^+-2\alpha} = 
	C|x|^{-p\beta^+}\le 	u^p -\rho |x|^{-p\beta^+} =U^p-\rho |x|^{-p\beta^+}\le -\mathcal{F}(U)$,
\end{center} 
up to changing $W$ by $c_0\,W$ for a suitable $c_0>0$.\smallskip

Hence $u(x)\geq W (x)$ by the Comparison Principle, as in the end of the proof of Lemma \ref{lemma varphi>0}, whenever $|x|> s_1R$, by taking $t_1\rightarrow \infty$. This gives us \eqref{eq inf log}.
The case $\beta^+\in (0,1)$ is analogous, see the proof of Theorem 1.3 in \cite{FQ}.

\smallskip

To finish the proof at the critical positive $p$ we use \eqref{eq inf log} and the rescaling $u_\sigma (x):=\sigma^{\frac{2\alpha}{p-1}}u(\sigma x)=\sigma^{\beta^+}u(\sigma x)$, to find that for $\sigma \geq s_1 R_1$,
\begin{equation*}
	u_\sigma (x) \geq C \log (1+\sigma) \quad \mbox{for} \;  x\in (B_{t_1/s_1}\setminus B_{1})\cap \mathcal{C}_{\omega_0}.
\end{equation*}
Thus for all $\sigma \geq s_1 R_1$,
\begin{equation*}
	-\mathcal{F}( u_\sigma) \geq  u_{\sigma}^p \geq C \log (1+\sigma)^{p-1} u_\sigma \quad \mbox{in} \; (B_{t_1/s_1}\setminus B_{1})\cap \mathcal{C}_{\omega_0}.
\end{equation*}
Therefore, by the definition of the first eigenvalue of $\mathcal{F}$ in $(B_{t/s}\setminus B_{1})\cap \mathcal{C}_{\omega_0}$, it is bounded below by $C (\log \sigma)^{p-1}$, a contradiction when one passes to the limit $\sigma \to+ \infty$, see \cite{Revista, QSX}.
\end{proof}

\begin{proof}[Proof of Theorem \ref{tm 0.3}] 
Let us first consider $p\in (\frac{\beta^-+2\alpha}{\beta^-},0]$. Note that Lemma \ref{lemaQ1}(ii) and \eqref{est D} implies 
\begin{equation*}
m(sr,tr,\omega_2)\leq\, \Phi^{-}_{s,t,\omega_2}(r)  \sup_{sr\le |x|\le tr, \,x\in \mathcal{C}_\omega} \phi^- \,\leq\, C r^{-\beta^-}.
\end{equation*}
Then Lemma~\ref{lemaQ1} yields
	\begin{center}
		$cr^{\frac{2\alpha}{1-p}}\le  m(sr,tr,\omega_2)\le C r^{ -\beta^-} $
	\end{center}
for large values of $r$, and so, since $ 2\alpha +\beta^- - p\beta^- >0$,
 \begin{center}
 	$0< c_0\le r^{- \frac{2\alpha + (1-p)\beta^-}{1-p} }\to 0$ \quad  as $r\to +\infty$,
 \end{center}
which is impossible. This concludes the proof in the subcritical nonnegative case.

\smallskip

As a second step, let $p=\frac{\beta^-+2\alpha}{\beta^-}$. We focus our attention on the case of the operator $\mathcal F= \M^-$.

We first note that if $u$ is a positive supersolution of \eqref{eq 0.002M-}, then we also have a solution of $\mathcal{M}^-(u)+u^p =0$ in $\mathcal{C}_\omega$, $u\ge 0$ in $\rn$, since $g(u)=u^p$ is a nonincreasing function of $u$ when $p<0$. W.l.g.\ we may assume from the beginning that $u$ is a positive solution which satisfies the lower bound $u(x)\ge C |x|^{-\beta^-}$ at infinity, by Lemma \ref{lemaQ1}(ii) and $\frac{2\alpha}{p-1}=\beta^-$.

On the other hand, by the Harnack inequality (see \cite[Theorem 11.1]{CS1}) we infer that
\begin{center}
	$M_1(r)\le C m_1(r)$ for large $r$.
\end{center}
Here,
	$m_1(r)=m(r/2,4r,\omega_1)=\inf_{1/2\le |x|\le 4,\, x\in \mathcal{C}_{\omega_1}} u_r (x)$, and  $M_1(r)=\sup_{1/2\le |x|\le 4,\, x\in \mathcal{C}_{\omega_1}} u_r (x)$ for some $\omega_1\subset\subset\omega$; recall that $u_r(x)=u(rx)$.

Next we consider, as in \cite{AScpde}, the function $w_r(x)=u_r(x)-\hat\varphi (r)\Psi (rx)$, where $\hat{\varphi}$ and $\Psi$ come from Lemma \ref{lemma rho bounded} with $s=1/2$, $t=4$. From that proof follows $w_r\ge 0$ in $B_4\backslash B_{1/2}$ for large $r$. Then $w_r$ satisfies
\begin{center}
$-\mathcal{M}^-(w_r)\ge u_r^{p}\ge c r^{-\beta^-} \chi_{\Omega_1} (x)$ \; in $\mathcal{C}_\omega$,
\end{center} where
$\Omega_1=(B_{4r}\backslash B_{\frac{r}{2}})\cap \mathcal{C}_{\omega_1}$. Now, from the definition of viscosity supersolution of $w(x):=w_r(x/r)$ we have, exactly as in the proof of Lemma \ref{lemaQ1}, for $\Omega_2=(B_{2r}\backslash B_{r})\cap \mathcal{C}_{\omega_2}$, that 
\begin{center}
	$cr^{-\beta^--2\alpha} \le Cr^{-2\alpha}\inf_{\Omega_2} w $, \;  i.e.\ \;
$\inf_{\Omega_2} w_r\ge cr^{-\beta^-}> \delta \phi^- (rx)$.
\end{center}

Finally, the function $\tilde w_r(x)=v_r(x)-\delta \phi^-  (rx)$ is a supersolution of $-\M^-(\tilde w_r)\ge 0$ in $\mathcal C_\omega$ by satisfying analogous properties as the function $w$ considered in Proposition \ref{propPhi2leC}, with $\Phi_0^-,\omega_0$ replaced by $1, \omega_2$ respectively.
So, following that proof we obtain $\tilde w_r\ge 0$ in $(B_{2r}\backslash B_{r})\cap \mathcal{C}_{\omega}$, up to taking a smaller $\delta$ if necessary. We then conclude that our initial $u_r$ satisfies 
$u_r \ge (\hat\varphi (r) +\delta) \Psi (rx) $ in $(B_{2r}\backslash B_{r})\cap \mathcal{C}_{\omega}$, from which follows $\hat\varphi (r/2) \ge \hat{\varphi} (r) +\delta $. By iterating we get $\hat{\varphi} (r)\to -\infty $ as $r\to \infty$, which is impossible. The proof of the critical negative case is then accomplished.

\smallskip

Finally, when $p<\frac{\beta^-+2\alpha}{\beta^-}$ the asymptotic behavior regarding the blow-up of supersolutions of \eqref{eq 0.002} follows from Lemma \ref{lemaQ1}(ii).
\end{proof}

  \bigskip
  
  \noindent {\bf Acknowledgements.}
  We would like to warmly thank P. Felmer for several interesting discussions on the  topic of this paper. 
We extend our sincere gratitude to the referee for his/her generous time and meticulous attention in the revision of this manuscript, the valuable  suggestions and for kindly bringing references \cite{FRbook, RSpotential} to our attention.

\smallskip
  
D.\ dos Prazeres was partially supported by CNPq grant 305680/2022-6 and Capes-Fapitec grant 88887.157906/2017-00.
  
   G.\ Nornberg was supported by Centro de Modelamiento Matemático FB210005, BASAL funds for centers of excellence from ANID-Chile, CMM-DIM, CNRS IRL 2807; by ANID Fondecyt grant 1220776, and by Vicerrectoría de Investigación y Desarrollo UI-001/21, U. de Chile.

A.\ Quaas was partially supported by ANID Fondecyt grant 1190282.

\end{document}